\documentclass{elsart}
\usepackage{graphics}
\usepackage{graphicx}
\usepackage{epsfig}
\usepackage{amssymb,amsmath}
\usepackage[colorlinks]{hyperref}
\newtheorem{lemma}{Lemma}[section]

\newtheorem{theorem}{Theorem}[section]
\newtheorem{proposition}{Proposition}[section]
\newtheorem{definition}{Definition}[section]
\newtheorem{remark}{Remark}[section]

\begin{document}
\linespread{1}
\begin{frontmatter}
%%%%%%%%%%%%%%%%%%%%%%%%%%%%%%%%%%%%%%%%%%%%%%%%%%%%%%%%%%%%%%%%%%%%%%%%%%%%%%%%%%%%%%%%%%%%%%%%%%%%%%%%%%%%%%%%%%%%%%%%%%%%%%%%%%%%%%%%%%%%%%%%%%%%%%%%%%%%%%%%%%%
\title{A mixed multifractal formalism for finitely many non Gibbs Frostman-like measures}
%%%%%%%%%%%%%%%%%%%%%%%%%%%%%%%%%%%%%%%%%%%%%%%%%%%%%%%%%%%%%%%%%%%%%%%%%%%%%%%%%%%%%%%%%%%%%%%%%%%%%%%%%%%%%%%%%%%%%%%%%%%%%%%%%%%%%%%%%%%%%%%%
%%%%%%%%%%%%%%%%
\author{Mohamed Menceur}
\address{Algerba and Number Theory Laboratory, Faculty of Mathematics, University of Sciences and Technology Houari Boumediene, BP 32 EL Alia 16111 Bab Ezzouar, Algiers, Algeria.}
\ead{m.m.m@live.fr}
%%%%%%%%%%%%%%%%
\author{Anouar Ben Mabrouk\corauthref{cor1}}
%\author{Anouar Ben Mabrouk\corauthref{cor1}\thanksref{label1}}
\address{Department of Mathematics, Higher Institute of Applied Mathematics and Informatics, Street of Assad Ibn Alfourat, 3100 Kairouan, Tunisia.}
%\address{{Research Unit of Algebra, Number Theory and Nonlinear Analysis, UR11ES50, Department of Mathematics, Faculty of Sciences, 5000 Monastir, Tunisia.}}
\ead{anouar.benmabrouk@fsm.rnu.tn}
%\thanks[label1]{{Department of Mathematics, Higher Institute of Applied Mathematics and Informatics, Street of Assad Ibn Alfourat, 3100 Kairouan, Tunisia.}}
\corauth[cor1]{Corresponding author.}
%%%%%%%%%%%%%%%%%%%%%%%%%%%%%%%%%%%%%%%%%%%%%%%%%%%%%%%%%%%%%%%%%%%%%%%%%%%%%%%%%%%%%%%%%%%%%%%%%%%%%%%%%%%%%%%%%%%%%%%%%%%%%%%%%%%%%%%%%%%%%%%%
%%%%%%%%%%%%%%%%%%%%%%%%%%%%%%%%%%%%%%%%%%%%%%%%%%%%%%%%%%%%%%%%%%%%%%%%%%%%%%%%%%%%%%%%%%%%%%%%%%%%%%%%%%%%%%%%%%%%%%%%%%%%%%%%%%%%%%%%%%%%%%%%
%%%%%%%%%%%%%%%%%%%%%%%%%%%%%%%%%%%%%%%%%%%%%%%%%%%%%%%%%%%%%%%%%%%%%%%%%%%%%%%%%%%%%%%%%%%%%%%%%%%%%%%%%%%%%%%%%%%%%%%%%%%%%%%%%%%%%%%%%%%%%%%%
%%%%%%%%%%%%%%%%%%%%%%%%%%%%%%%%%%%%%%%%%%%%%%%%%%%%%%%%%%%%%%%%%%%%%%%%%%%%%%%%%%%%%%%%%%%%%%%%%%%%%%%%%%%%%%%%%%%%%%%%%%%%%%%%%%%%%%%%%%%%%%%%
\begin{abstract}
The multifractal formalism for measures hold whenever the existence of corresponding Gibbs-like measures supported on the singularities sets holds. In the present work we tried to relax such a hypothesis and introduce a more general framework of mixed (and thus single) multifractal analysis where the measures constructed on the singularities sets are not Gibbs but controlled by an extra-function allowing the multifractal formalism to hold. We fall on the classical case by a particular choice of such afunction.
\end{abstract}
\end{frontmatter}
\tableofcontents
\section{Introduction}
The multifractal analysis of a single measure passes through its local dimension or its H\"older exponent. For a measure $\mu$ eventually Borel and finite on $\mathbb{R}^d$ and $x\in\,\hbox{supp}(\mu)$, the local dimension of $\mu$ at the point $x$ is defined by
$$
{\alpha}_{\mu}(x)=\displaystyle\lim_{r\downarrow0}\displaystyle\frac{\log(\mu(B(x,r)))}{\log\,r}
$$
when such a limit exists. The next step concerns the geometric study of the $\alpha$-singularity set of the measure $\mu$ defined by
$$
X(\alpha)=\{\,x\in\hbox{supp}(\mu)\,;\,\,\alpha_{\mu}(x)=\alpha\,\}
$$
by means of its Hausdorff dimension
$$
d(\alpha)=\,\hbox{dim}\,X(\alpha)
$$
which defines the so-called spectrum of singularities. This means that the study of the behaviour of the measure is transformed into a study of sets where the focuses may somehow forget about the measure and its point-wise character and falls in set theory and the suitable coverings that permits the computation of the Hausdorff dimension. For a subset $E\subset\mathbb{R}^d$ and $\alpha\geq0$, the $\alpha$-Hausdorff measure is defined by
$$
\mathcal{H}^\alpha(E)=\displaystyle\lim_{\varepsilon\downarrow0}\left(\inf\displaystyle\sum_j(diam(U_j))^\alpha\,\right)
$$
where the inf is taken over all coverings of $E$ with subsets $U_j$, $j\in\mathbb{N}$ such that $diam(U_j)\leq\varepsilon$.

However, some geometric sets are essentially known by means of measures that are supported by them, i.e., given a set $E$ and a measure $\mu$, the quantity $\mu(E)$ may be computed as the maximum value $\mu(F)$ for all subsets $F\subset E$. So, contrarily to the previous idea, we mathematically forget the geometric set structure of $E$ and focus instead on the properties of the measure $\mu$. The set $E$ is thus partitioned into $\alpha$-level sets relatively to the regularity exponent of $\mu$ into subsets $X(\alpha)$.

This makes the including of the measure $\mu$ into the computation of the Hausdorff (or fractal) dimension and thus into the definition of the Hausdorff measure a necessity to understand more the geometry of the set simultaneously with the behaviour of the measure that is supported on. One step ahead in this direction has been conducted by Olsen in \cite{olsen1} where the author introduced multifractal generalisations of the fractal dimensins such as Hausdorff, packing and Bouligand ones by considering general variants of measures. For a Borel probability measure $\mu$ on $\mathbb{R}^d$, a nonempty set $E\subseteq\mathbb{R}^d$ and $q,\,t\in\mathbb{R}$, he considered the pre-mesure
$$
\overline{\mathcal{H}}_\mu^{q,t}(E)=\displaystyle\lim_{\epsilon\downarrow0}\bigl(\displaystyle\inf\{\,\sum_i(\mu(B(x_i,r_i)))^q(2r_i)^t\,\}\bigr),
$$
where the inf is taken over the set of all centered $\epsilon$-coverings of $E$, and for the empty set, $\overline{\mathcal{H}}_{\mu,\epsilon}^{q,t}(\emptyset)=0$. This yields next the measure
$$
\mathcal{H}_\mu^{q,t}(E)=\displaystyle\sup_{F\subseteq\,E}\overline{\mathcal{H}}_\mu^{q,t}(F).
$$
Similarly, the following pre-measure is considered.
$$
\overline{\mathcal{P}}_\mu^{q,t}(E)=\displaystyle\lim_{\epsilon\downarrow0}\bigl(\displaystyle\sup\{\,\sum_i(\mu(B(x_i,r_i)))^q(2r_i)^t\,\}\bigr),
$$
where the sup is taken over the set of all centered $\epsilon$-packings of $E$. For the empty set, we set as usual $\overline{\mathcal{P}}_{\mu,\epsilon}^{q,t}(\emptyset)=0$. This yields next the measure
$$
\mathcal{P}_\mu^{q,t}(E)=\displaystyle\inf_{E\subseteq\,\cup_iE_i}\sum_i\overline{\mathcal{P}}_\mu^{q,t}(E_i).
$$
In \cite{olsen1}, it has been proved that the measures $\mathcal{H}_{\mu}^{q,t}$, $\mathcal{P}_\mu^{q,t}$ and the pre-measure $\overline{\mathcal{P}}_\mu^{q,t}$ assign in a usual way a dimension to every set $E\subseteq\mathbb{R}^d$ called respectively multifractal generalizations of the Hausdorff dimension ($\hbox{dim}_\mu^q(E)$,), the packing dimension ($\hbox{Dim}_\mu^q(E)$) and the logarithmic index ($\Delta_\mu^q(E)$) of the set $E$. These quantities satisfies the cut-off relations
$$
\mathcal{H}_\mu^{q,t}(E)=+\infty\;\hbox{for}\;t<\hbox{dim}_\mu^q(E)\;\;\hbox{and}\;\;0\;\hbox{for}\;t>\hbox{dim}_\mu^q(E).
$$
$$
\mathcal{P}_\mu^{q,t}(E)=+\infty\;\hbox{for}\;t<\hbox{Dim}_\mu^q(E)\;\;\hbox{and}\;\;0\;\hbox{for}\;t>\hbox{Dim}_\mu^q(E).
$$
$$
\overline{\mathcal{P}}_\mu^{q,t}(E)=+\infty\;\hbox{for}\;t<\Delta_\mu^q(E)\;\;\hbox{and}\;\;0\;\hbox{for}\;t>\Delta_\mu^q(E).
$$
In \cite{olsen1}, the characteristics of these functions such as monotony, convexity, lower and upper bounds have been studied. Next, to come back to the essential problem in multifractal formalism which consists in the computation of the spectrum of singularities $d(\alpha)$, the author proved that such generalizations may lead to a multifractal formalism but when some bad restriction for the single measure $\mu$ toke place. By assuming that $\mu$ belongs to the whole class of Gibbs-like measures the multifractal formalism has been proved to hold. This was one motivation that led us to develop the present paper where such a restriction has been avoided.

Next, a first step in the direction of the mixed multifractal analysis the same author already affected by the restriction on the single measure $\mu$ developed a mixed multifractal analysis for one very restrictive class of measures known as the self affine measures \cite{olsen3} dealing precisely with R\'enyi dimensions for finitely many self similar measures. It was one step ahead but in a restrictive case. This study has been the motivation of our paper \cite{Bmabrouk3} where we developed a general mixed analysis for vector valued measures by proving some results for general measures and some ones for special classes. However, we noticed that the hypothesis of Gibbs-like measures is somehow not possible to avoid and thus by contouring such hypothesis with some extra-hypothesis on the measure and by proving a general mixed large deviation formalism a mixed multifractal formalism has been proved there also. This study itself has been one motivation behind the one developed in \cite{Bmabrouk4} and \cite{Bmabrouk5} where a mixed multifractal analysis inspired from the one for measures has been developed in the functional case. By concentring a vector valued Gibbs-like measure on the singularities set of finitely and simultaneously many functions, we introduced a mixed multifractal formalism for functions. General results for almost all functions have been proved and a mixed multifractal formalism have been proved for self similar quasi self similar functions as well as their superpositions (which are not self similar neither quasi self similar). For more details and backgrounds on multifractal analysis as well as the mixed generalizations the readers may be referred also to the following essential references \cite{Bmabrouk1}, \cite{Bmabrouk2}, \cite{olsen2b}, \cite{Xu-Xu}, \cite{Xu-Xu-Zhong}, \cite{Ye}, \cite{Yuan}, \cite{Zeng-Yuan-Xui}, \cite{Zhou-Feng}, \cite{Zhu-Zhou}.

In the present paper we are concerned with the introduction of a multifractal analysis in a mixed case (but which can already adapted to single cases) where the hypothesis of the existence of Gibbs-like and/or doubling measures supported by the singularities sets is relaxed. We aim to consider some cases of simultaneous behaviors of measures where the local H\"older behaviour is controlled by special and suitable function that allow the extra-hypothesis of Gibbs-like measures not to be necessary.

The present work will be organized as follows. The next section concerns the introduction of the new variant of the mixed multifractal generalizations of Hausdorff and packing measures and dimensions relatively to the control function $\varphi$. Section 3 is devoted to the mixed multifractal generalization of Bouligand-Minkowsky or R\'enyi dimension already with the control $\varphi$. In section 4, a mixed multifractal formalism associated to the mixed multifractal generalizations of Hausdorff and packing measures and dimensions introduced in section 2 is proved in some case based on a generalization of the well known large deviation formalism and where no extra-hypothesis of Gibbs-like measures existence is assumed.
\section{$\varphi$-mixed multifractal generalizations of Hausdorff and packing measures and dimensions}
The purpose of this section is to present our ideas about mixed multifractal generalizations of Hausdorff and packing measures and dimensions. Consider a vector valued measure $\mu=(\mu_1,\,\mu_2,\dots,\mu_k)$ composed of probability measures on $\mathbb{R}^d$. We aim to study the simultaneous scaling behavior of $\mu$ relatively to an exponential density function. Let $\varphi:\mathbb{R}_+\rightarrow\mathbb{R}$ be such that
\begin{equation}\label{assumptiononphi}
	\varphi\;\hbox{is non-decreasing and}\;\;\varphi(r)<0\;\hbox{for}\;r\;\hbox{small enough}.
\end{equation}
The mixed generalized multifractal Hausdorff $\varphi$-measure is defined as follows. For $x\in\mathbb{R}^d$ and $r>0$ we denote $B(x,r)$ the ball of radius $r$ and center $x$. We denote next
$$
\mu(B(x,r))\equiv\bigl(\mu_1(B(x,r)),\dots,\,\mu_k(B(x,r))\bigr)
$$
and the product
$$
(\mu(B(x,r)))^q\equiv(\mu_1(B(x,r)))^{q_1}\dots(\mu_k(B(x,r)))^{q_k}.
$$
Let $E\subseteq\mathbb{R}^d$ be a nonempty set and $\epsilon>0$. Let also $q=(q_1,q_2,\dots,q_k)\in\mathbb{R}^k$ and $t\in\mathbb{R}$ and consider the quantity
$$
\overline{\mathcal{H}}_{\mu,\varphi,\epsilon}^{q,t}(E)=\displaystyle\inf\{\,\sum_i(\mu(B(x_i,r_i)))^qe^{t\varphi(r_i)}\,\},
$$
where the inf is taken over the set of all centered $\epsilon$-coverings of $E$, and for the empty set, $\overline{\mathcal{H}}_{\mu,\epsilon}^{q,t}(\emptyset)=0$. It consists of a non increasing function of the variable $\varepsilon$. We denote thus
$$
\overline{\mathcal{H}}_{\mu,\varphi}^{q,t}(E)=\displaystyle\lim_{\epsilon\downarrow0}\overline{\mathcal{H}}_{\mu,\varphi,\epsilon}^{q,t}(E)
=\sup_{\delta>0}\overline{\mathcal{H}}_{\mu,\varphi,\epsilon}^{q,t}(E).
$$
Let finally
$$
\mathcal{H}_{\mu,\varphi}^{q,t}(E)=\displaystyle\sup_{F\subseteq\,E}\overline{\mathcal{H}}_{\mu,\varphi}^{q,t}(F).
$$
\begin{lemma}\label{anouarlemme1}
	$\mathcal{H}_{\mu,\varphi}^{q,t}$ is an outer metric measure on $\mathbb{R}^d$.
\end{lemma}
\hskip-10pt{\bf Proof.} We will prove firstly that $\mathcal{H}_{\mu,\varphi}^{q,t}$ is an outer measure. This means that
\begin{description}
	\item [i.] $\mathcal{H}_{\mu,\varphi}^{q,t}(\emptyset)=0$.
	\item[ii.] $\mathcal{H}_{\mu,\varphi}^{q,t}$ is monotone, i.e.
	$$
	\mathcal{H}_{\mu,\varphi}^{q,t}(E)\leq\mathcal{H}_{\mu,\varphi}^{q,t}(F),
	$$
	whenever $E\subseteq\,F\subseteq\mathbb{R}^d$.
	\item[iii.] $\mathcal{H}_{\mu,\varphi}^{q,t}$ is sub-additive, i.e.
	$$
	\mathcal{H}_{\mu,\varphi}^{q,t}(\displaystyle\bigcup_nA_n)\leq\displaystyle\sum_n\mathcal{H}_{\mu,\varphi}^{q,t}(A_n).
	$$
\end{description}
The first item is obvious. Let us prove (ii). Let $E\subseteq\,F$ be nonempty subsets of $\mathbb{R}^d$. We have
$$
\mathcal{H}_{\mu,\varphi}^{q,t}(E)=\displaystyle\sup_{A\subseteq\,E}\overline{\mathcal{H}}_{\mu,\varphi}^{q,t}(A)
\leq\displaystyle\sup_{A\subseteq\,F}\overline{\mathcal{H}}_{\mu,\varphi}^{q,t}(A)=\mathcal{H}_{\mu,\varphi}^{q,t}(F).
$$
We next prove (iii). If the right hand term is infinite, the inequality is obvious. So, assume that it is finite. Let $(E_n)_n$ be a countable family of subsets $E_i\subseteq\mathbb{R}^d$ for which $\displaystyle\sum_n\mathcal{H}_{\mu,\varphi}^{q,t}(E_n)<\infty$. Let also $\epsilon,\delta>0$ and $(B(x_{ni},r_{ni}))_i$ a centered $\epsilon$-covering of $E_n$ satisfying
$$
\displaystyle\sum_i\mu(B(x_{ni},r_{ni}))^qe^{t\varphi(r_{ni})}\leq\overline{\mathcal{H}}_{\mu,\varphi,\varepsilon}^{q,t}(E_n)+\frac{\delta}{2^n}.
$$
The whole set $(B(x_{ni},r_{ni}))_{n,i}$ is a centered $\epsilon$-covering of the whole union $\displaystyle\bigcup_nE_n$. As a consequence,
$$
\begin{array}{lll}\displaystyle\overline{\mathcal{H}}_{\mu,\varphi,\epsilon}^{q,t}(\bigcup_nE_n)
&\leq&\displaystyle\sum_n\sum_i(\mu(B(x_{ni},r_{ni})))^qe^{t\varphi(r_{ni})}\hfill\cr\medskip
&\leq&\displaystyle\sum_n\bigl(\overline{\mathcal{H}}_{\mu,\varphi,\epsilon}^{q,t}(E_n)+\frac{\delta}{2^n}\bigr)\hfill\cr\medskip
&\leq&\displaystyle\sum_n\bigl(\overline{\mathcal{H}}_{\mu,\varphi}^{q,t}(E_n)+\frac{\delta}{2^n}\bigr)\hfill\cr\medskip
&\leq&\displaystyle\sum_n\mathcal{H}_{\mu,\varphi}^{q,t}(E_n)+\delta.\end{array}$$
Having $\epsilon$ and $\delta$ going towards 0, we obtain
$$
\overline{\mathcal{H}}_{\mu,\varphi}^{q,t}(\displaystyle\bigcup_nE_n)\leq\displaystyle\sum_n\mathcal{H}_{\mu,\varphi}^{q,t}(E_n).
$$
Let next a set $F$ covered with the countable set $(A_n)_n$. That is $F\subseteq\displaystyle\bigcup_nA_n$. We have
$$
\begin{array}{lll}\overline{\mathcal{H}}_{\mu,\varphi}^{q,t}(F)
&=&\overline{\mathcal{H}}_{\mu,\varphi}^{q,t}\biggl(\displaystyle\bigcup_n(A_n\cap\,F)\biggr)\hfill\cr\medskip
&\leq&\displaystyle\sum_n\mathcal{H}_{\mu,\varphi}^{q,t}(A_n\cap\,F)\hfill\cr\medskip
&\leq&\displaystyle\sum_n\mathcal{H}_{\mu,\varphi}^{q,t}(A_n).\end{array}
$$
Taking the sup on $F$, we obtain
$$
\mathcal{H}_{\mu,\varphi}^{q,t}(F)\leq\displaystyle\sum_n\mathcal{H}_{\mu,\varphi}^{q,t}(A_n).
$$
We now prove that $\mathcal{H}_{\mu,\varphi}^{q,t}$ is metric. Let $A,B$ subsets of $\mathbb{R}^d$. We recall that the Hausdorff distance $d(A,B)$ is defined by
$$
d(A,B)=inf\{|x-y|;\;\;x\in\,A\,\;\;y\in\,B\}.
$$
Assume so that  $d(A,B)>0$, and that
$$
\mathcal{H}_{\mu,\varphi}^{q,t}(A\cup\,B)<\infty.
$$
Let next $0<\delta<d(A,B)$, $\varepsilon>0$, $F_1\subseteq\,A$, $F_2\subseteq\,B$ and $(B(x_i,r_i))_i$ a centered $\delta$-covering of the set $F_1\cup\,F_2$ and such that
$$
\overline{\mathcal{H}}_{\mu,\varphi,\delta}^{q,t}(F_1\cup\,F_2)\leq\displaystyle\sum_i(\mu(B(x_i,r_i)))^qe^{t\varphi(r_{i})}
\leq\overline{\mathcal{H}}_{\mu,\varphi,\delta}^{q,t}(F_1\cup\,F_2)+\varepsilon.
$$
This is always possible from the definition of $\overline{\mathcal{H}}_{\mu,\varphi,\delta}^{q,t}(F_1\cup\,F_2)$.
Denote next the index sets
$$
I=\{\,i;\;\;B(x_i,r_i)\cap\,F_1\not=\emptyset\,\}\quad\hbox{and}\quad\,J=\{\,i;\;\;B(x_i,r_i)\cap\,F_2\not=\emptyset\,\}.
$$
The countable sets $(B(x_i,r_i))_{i\in\,I}$ and $(B(x_i,r_i))_{i\in\,J}$ are centered $\delta$-coverings of $F_1$ and $F_2$ respectively. Consequently,
$$
\begin{array}{lll}\overline{\mathcal{H}}_{\mu,\varphi,\delta}^{q,t}(F_1)+\overline{\mathcal{H}}_{\mu,\varphi,\delta}^{q,t}(F_2)
&\leq&\displaystyle\sum_{i\in\,I}(\mu(B(x_i,r_i)))^qe^{t\varphi(r_{i})}\hfill\cr\medskip
&&+\displaystyle\sum_{i\in\,J}\mu(B(x_i,r_i))^qe^{t\varphi(r_{i})}\hfill\cr\medskip
&=&\displaystyle\sum_i(\mu(B(x_i,r_i)))^q(2r_i)^t\hfill\cr\medskip
&\leq&\overline{\mathcal{H}}_{\mu,\varphi,\delta}^{q,t}(F_1\cup\,F_2)+\varepsilon.\end{array}
$$
As a result,
$$
\overline{\mathcal{H}}_{\mu,\varphi}^{q,t}(F_1)+\overline{\mathcal{H}}_{\mu,\varphi}^{q,t}(F_2)
\leq\overline{\mathcal{H}}_{\mu,\varphi}^{q,t}(F_1\cup\,F_2)
+\varepsilon\leq\mathcal{H}_{\mu,\varphi}^{q,t}(A\cup\,B)+\varepsilon.
$$
When $\varepsilon\downarrow0$ and taking the sup on the sets $F_1\subseteq\,A$ and $F_2\subseteq\,B$, we obtain
$$
\mathcal{H}_{\mu,\varphi}^{q,t}(A\cup\,B)\geq\mathcal{H}_{\mu,\varphi}^{q,t}(A)+\mathcal{H}_{\mu,\varphi}^{q,t}(B).
$$
The inequality
$$
\mathcal{H}_{\mu,\varphi}^{q,t}(A\cup\,B)\leq\mathcal{H}_{\mu,\varphi}^{q,t}(A)+\mathcal{H}_{\mu,\varphi}^{q,t}(B).
$$
follows from the sub-additivity property of the measure $\mathcal{H}_{\mu,\varphi}^{q,t}$.
\begin{lemma}\label{Hmuphi-1}
	Let $(A_n)_n$ be a non-decreasing sequence of subsets in $\mathbb{R}^d$ and denote $A=\displaystyle\bigcup_nA_n$. Assume further that $d(A_n,A\setminus A_{n+1})>0$ for all $n$. Then,
	$$
	\mathcal{H}_{\mu,\varphi}^{q,t}(A)=\displaystyle\lim_{n\rightarrow+\infty}\mathcal{H}_{\mu,\varphi}^{q,t}(A_n).
	$$
\end{lemma}
\hskip-10pt\textbf{Proof.} The result is obvious if the limit above is infinite. So assume that it is finite and denote for $k\in\mathbb{N}$, $C_k=A_{k+1}\setminus A_k$. We then observe that
$$
d(C_j,C_p)>0,\;\forall j,p\,;\;|j-p|\geq2\quad\hbox{and}\quad\,A=A_n\cup\bigl(\displaystyle\bigcup_{k\geq n+1}C_k\bigr),\;\forall\,n.
$$
Therefore,
$$
\mathcal{H}_{\mu,\varphi}^{q,t}(A)\leq
\mathcal{H}_{\mu,\varphi}^{q,t}(A_n)+\underbrace{\displaystyle\sum_{k\geq n+1}\mathcal{H}_{\mu,\varphi}^{q,t}(C_k).}_{R_n}
$$
Now, it is straightforward that for all $n$ we have
$$
\displaystyle\sum_{k=0}^n\mathcal{H}_{\mu,\varphi}^{q,t}(C_{2k})=\mathcal{H}_{\mu,\varphi}^{q,t}\left(\bigcup_{k=0}^nC_{2k}\right)
\leq\mathcal{H}_{\mu,\varphi}^{q,t}(A_{2n+2})\leq\mathcal{H}_{\mu,\varphi}^{q,t}(A)<\infty
$$
and similarly,
$$
\displaystyle\sum_{k=0}^n\mathcal{H}_{\mu,\varphi}^{q,t}(C_{2k+1})=\mathcal{H}_{\mu,\varphi}^{q,t}\left(\bigcup_{k=0}^nC_{2k+1}\right)
\leq\mathcal{H}_{\mu,\varphi}^{q,t}(A_{2n+1})\leq\mathcal{H}_{\mu,\varphi}^{q,t}(A)<\infty.
$$
Hence, $R_n$ is the rest of a convergent series, so it goes to 0 as $n$ goes to infinity. Consequently,
$$
\mathcal{H}_{\mu,\varphi}^{q,t}(A)\leq\displaystyle\lim_{n\rightarrow+\infty}\mathcal{H}_{\mu,\varphi}^{q,t}(A_n).
$$
The reciprocal inequality is obvious.
\begin{lemma}
	Borel sets are $\mathcal{H}_{\mu,\varphi}^{q,t}$-measurable.
\end{lemma}
\hskip-10pt\textbf{Proof.} Let $B$ be a Borel subset of $\mathbb{R}^d$ (a closed subset for example), $E\subset\mathbb{R}^d$ and denote for $n\in\mathbb{N}$,
$$
B_n=\{x\in E\,;\;d(x,B)\geq\frac{1}{n}\}.
$$
It consists of a non-decreasing sequence of subsets of $\mathbb{R}^d$ satisfying further that
$$
E\setminus B=\bigcup_nB_n\quad\hbox{and}\quad\,d(B_n,(E\setminus B)\setminus B_{n+1})>0,\;\forall\,n.
$$
Hence, Lemma \ref{Hmuphi-1} yields that
$$
\mathcal{H}_{\mu,\varphi}^{q,t}(E\setminus B)=\displaystyle\lim_{n\rightarrow+\infty}\mathcal{H}_{\mu,\varphi}^{q,t}(B_n).
$$
Observe now that
$$
\mathcal{H}_{\mu,\varphi}^{q,t}(E)\geq\mathcal{H}_{\mu,\varphi}^{q,t}(E\cap B)+\mathcal{H}_{\mu,\varphi}^{q,t}(B_n),\;\forall\,n.
$$
When $n\rightarrow\infty$, we get
$$
\mathcal{H}_{\mu,\varphi}^{q,t}(E)\geq\mathcal{H}_{\mu,\varphi}^{q,t}(E\cap B)+\mathcal{H}_{\mu,\varphi}^{q,t}(E\setminus B).
$$
\begin{definition}
	The restriction of $\mathcal{H}_{\mu,\varphi}^{q,t}$ on Borel sets is called the mixed generalized Hausdorff measure on $\mathbb{R}^d$.
\end{definition}
Now, we define the mixed generalized multifractal packing measure. We use already the same notations as previously. Let
$$
\overline{\mathcal{P}}_{\mu,\varphi,\epsilon}^{q,t}(E)=\displaystyle\sup\{\,\sum_i(\mu(B(x_i,r_i)))^qe^{t\varphi(r_{i})}\,\}
$$
where the sup is taken over the set of all centered $\epsilon$-packings of $E$. For the empty set, we set as usual $\overline{\mathcal{P}}_{\mu,\varphi,\epsilon}^{q,t}(\emptyset)=0$. Next, we consider the limit as $\epsilon\downarrow0$,
$$
\overline{\mathcal{P}}_{\mu,\varphi}^{q,t}(E)=\displaystyle\lim_{\epsilon\downarrow0}\overline{\mathcal{P}}_{\mu,\varphi,\epsilon}^{q,t}(E)
=\inf_{\delta>0}\overline{\mathcal{P}}_{\mu,\varphi,\epsilon}^{q,t}(E)
$$
and finally,
$$
\mathcal{P}_{\mu,\varphi}^{q,t}(E)=\displaystyle\inf_{E\subseteq\,\cup_iE_i}\sum_i\overline{\mathcal{P}}_{\mu,\varphi}^{q,t}(E_i).
$$
\begin{lemma}\label{anouarlemme2}
	$\mathcal{P}_{\mu,\varphi}^{q,t}$ is an outer metric measure on $\mathbb{R}^d$.
\end{lemma}
\hskip-10pt The proof of this lemma is more specific than Lemma \ref{anouarlemme1} and uses the following result.
\begin{equation}\label{anouarclaim1}
	\overline{\mathcal{P}}_{\mu,\varphi}^{q,t}(A\cup\,B)=\overline{\mathcal{P}}_{\mu,\varphi}^{q,t}(A)
	+\overline{\mathcal{P}}_{\mu,\varphi}^{q,t}(B),\;\;\hbox{whenever}\;d(A,B)>0.
\end{equation}
Indeed, let $0<\epsilon<\displaystyle\frac{1}{2}d(A,B)$ and $(B(x_i,r_i))_i$ be a centered $\epsilon$-packing of the union $A\cup\,B$. It can be divided into two parts $I$ and $J$,
$$
(B(x_i,r_i))_i=\Bigl(B(x_i,r_i)\Bigr)_{i\in\,I}\bigcup\Bigl(B(x_i,r_i)\Bigr)_{i\in\,J}
$$
where
$$
\forall\,i\in\,I,\;\;B(x_i,r_i)\cap\,B=\emptyset\quad\hbox{and}\quad\forall\,i\in\,J,\;\;B(x_i,r_i)\cap\,A=\emptyset.
$$
Therefore, $(B(x_i,r_i))_{i\in\,I}$ is a centered $\epsilon$-packing of $A$ and $(B(x_i,r_i))_{i\in\,J}$ is a centered $\epsilon$-packing of the union $B$. Hence,
$$
\begin{array}{lll}
\displaystyle\sum_i(\mu(B(x_i,r_i)))^qe^{t\varphi(r_{i})}
&=&\displaystyle\sum_{i\in\,I}(\mu(B(x_i,r_i)))^qe^{t\varphi(r_{i})}\hfill\cr\medskip
&&+\displaystyle\sum_{i\in\,I}(\mu(B(x_i,r_i)))^qe^{t\varphi(r_{i})}.\hfill
\end{array}
$$
Now, it is straightforward that
$$
\displaystyle\sum_{i\in\,I}(\mu(B(x_i,r_i)))^qe^{t\varphi(r_{i})}\leq\overline{\mathcal{P}}_{\mu,\varphi,\epsilon}^{q,t}(A)
$$
and
$$
\displaystyle\sum_{i\in\,I}(\mu(B(x_i,r_i)))^qe^{t\varphi(r_{i})}\leq\overline{\mathcal{P}}_{\mu,\varphi,\epsilon}^{q,t}(B).
$$
Consequently,
$$
\overline{\mathcal{P}}_{\mu,\varphi,\epsilon}^{q,t}(A\cup\,B)\leq\overline{\mathcal{P}}_{\mu,\varphi,\epsilon}^{q,t}(A)
+\overline{\mathcal{P}}_{\mu,\varphi,\epsilon}^{q,t}(B)
$$
and thus the limit for $\epsilon\downarrow0$ gives
$$
\overline{\mathcal{P}}_{\mu,\varphi}^{q,t}(A\cup\,B)\leq\overline{\mathcal{P}}_{\mu,\varphi}^{q,t}(A)+\overline{\mathcal{P}}_{\mu,\varphi}^{q,t}(B).
$$
The converse is more easier and it states that $\overline{\mathcal{P}}_{\mu,\varphi,\epsilon}^{q,t}$ and next $\overline{\mathcal{P}}_{\mu,\varphi}^{q,t}$ are sub-additive. Let $(B(x_i,r_i))_i$ be a centered $\epsilon$-packing of $A$ and $(B(y_i,r_i))_i$ be a centered $\epsilon$-packing of $B$. The union $\Bigl(B(x_i,r_i)\Bigr)_i\bigcup\Bigl(B(y_i,r_i)\Bigr)_i$ is a centered $\epsilon$-packing of $A\cup\,B$. So that
$$
\overline{\mathcal{P}}_{\mu,\varphi,\epsilon}^{q,t}(A\cup\,B)\geq\displaystyle\sum_i(\mu(B(x_i,r_i)))^qe^{t\varphi(r_{i})}
+\displaystyle\sum_i(\mu(B(y_i,r_i)))^qe^{t\varphi(r_{i})}.
$$
Taking the sup on $(B(x_i,r_i))_i$ as a centered $\epsilon$-packing of $A$ and next the sup on $(B(y_i,r_i))_i$ as a centered $\epsilon$-packing of $B$, we obtain
$$
\overline{\mathcal{P}}_{\mu,\varphi,\epsilon}^{q,t}(A\cup\,B)\geq\overline{\mathcal{P}}_{\mu,\varphi,\epsilon}^{q,t}(A)
+\overline{\mathcal{P}}_{\mu,\varphi,\epsilon}^{q,t}(B)
$$
and thus the limit for $\epsilon\downarrow0$ gives
$$
\overline{\mathcal{P}}_{\mu,\varphi}^{q,t}(A\cup\,B)\geq\overline{\mathcal{P}}_{\mu,\varphi}^{q,t}(A)+\overline{\mathcal{P}}_{\mu,\varphi}^{q,t}(B).
$$
\textbf{Proof of Lemma \ref{anouarlemme2}.} We shall prove as previously that
\begin{description}
	\item [i.] $\mathcal{P}_{\mu,\varphi}^{q,t}(\emptyset)=0$.
	\item[ii.] $\mathcal{P}_{\mu,\varphi}^{q,t}$ is monotone, i.e.
	$$
	\mathcal{P}_{\mu,\varphi}^{q,t}(E)\leq\mathcal{P}_{\mu,\varphi}^{q,t}(F),
	$$
	whenever $E\subseteq\,F\subseteq\mathbb{R}^d$.
	\item[iii.] $\mathcal{P}_{\mu,\varphi}^{q,t}$ is sub-additive, i.e.
	$$
	\mathcal{P}_{\mu,\varphi}^{q,t}(\displaystyle\bigcup_nA_n)\leq\displaystyle\sum_n\mathcal{P}_{\mu,\varphi}^{q,t}(A_n).
	$$
\end{description}
The first item is immediate from the definition of $\mathcal{P}_{\mu,\varphi}^{q,t}(\emptyset)=0$. Let $E\subseteq\,F$ be subsets of $\mathbb{R}^d$. We have
$$
\begin{array}{lll}
\medskip\mathcal{P}_{\mu,\varphi}^{q,t}(E)
&=&\displaystyle\inf_{E\subseteq\displaystyle\bigcup_iE_i}\displaystyle\sum_i\overline{\mathcal{P}}_{\mu,\varphi}^{q,t}(E_i)\hfill\cr\medskip
&\leq&\displaystyle\inf_{F\subseteq\displaystyle\bigcup_iE_i}\displaystyle\sum_i\overline{\mathcal{P}}_{\mu,\varphi}^{q,t}(E_i)\hfill\cr\medskip
&=&\mathcal{P}_{\mu,\varphi}^{q,t}(F).\hfill
\end{array}
$$
So is the item {\bf ii}. Let next $(A_n)_n$ a countable set of subsets of $\mathbb{R}^d$, $\varepsilon>0$ and for each $n$, $(E_{ni})_i$ be a covering of
$A_n$ such that
$$
\displaystyle\sum_i\overline{\mathcal{P}}_{\mu,\varphi}^{q,t}(E_{ni})\leq\mathcal{P}_{\mu,\varphi}^{q,t}(A_n)+\frac{\varepsilon}{2^n}.
$$
It follows for all $\varepsilon>0$ that
$$
\begin{array}{lll}
\medskip\mathcal{P}_{\mu,\varphi}^{q,t}(\displaystyle\bigcup_nA_n)
&\leq&\displaystyle\sum_n\displaystyle\sum_i\overline{\mathcal{P}}_{\mu,\varphi}^{q,t}(E_{ni})\hfill\cr\medskip
&\leq&\displaystyle\sum_n\mathcal{P}_{\mu,\varphi}^{q,t}(A_n)+\varepsilon.\hfill
\end{array}
$$
Hence,
$$
\mathcal{P}_{\mu,\varphi}^{q,t}(\displaystyle\bigcup_nA_n)\leq\displaystyle\sum_n\mathcal{P}_{\mu,\varphi}^{q,t}(A_n).
$$
So is the item {\bf iii}. We now prove that $\mathcal{P}_{\mu,\varphi}^{q,t}$ is metric. Let $A,B$ subsets of $\mathbb{R}^d$ be such that $d(A,B)>0$. We shall prove that
$$
\mathcal{P}_{\mu,\varphi}^{q,t}(A\cup\,B)=\mathcal{P}_{\mu,\varphi}^{q,t}(A)+\mathcal{P}_{\mu,\varphi}^{q,t}(B).
$$
Since $\mathcal{P}_{\mu,\varphi}^{q,t}$ is an outer measure, it suffices to show that
$$
\mathcal{P}_{\mu,\varphi}^{q,t}(A\cup\,B)\geq\mathcal{P}_{\mu,\varphi}^{q,t}(A)+\mathcal{P}_{\mu,\varphi}^{q,t}(B).
$$
Of course, if the left hand term is infinite, the inequality is obvious. So, suppose that it is finite. For $\varepsilon>0$, there exists a covering $(E_i)_i$ of the union set $A\cup\,B$ such that
$$
\displaystyle\sum_i\overline{\mathcal{P}}_{\mu,\varphi}^{q,t}(E_i)\leq\mathcal{P}_{\mu,\varphi}^{q,t}(A\cup\,B)+\varepsilon.
$$
By denoting $F_i=A\cap\,E_i$ and $H_i=B\cap\,E_i$, we get countable coverings $(F_i)_i$ of $A$ and $(H_i)_i$ for $B$ respectively. Furthermore, $F_i\cap\,H_j=\emptyset$ pour all $i$ and $j$. Consequently,
$$
\mathcal{P}_{\mu,\varphi}^{q,t}(A)+\mathcal{P}_{\mu,\varphi}^{q,t}(B)\leq\displaystyle\sum_i(\overline{\mathcal{P}}_{\mu,\varphi}^{q,t}(F_i)
+\overline{\mathcal{P}}_{\mu,\varphi}^{q,t}(H_i)).
$$
Since $d(A,B)>0$, $F_i\subset\,A$ and $H_i\subset\,B$, it follows that $d(F_i,H_j)>0$ for all $i$ and $j$. Hence, claim \ref{anouarclaim1} affirms that
$$
\overline{\mathcal{P}}_{\mu,\varphi}^{q,t}(E_i)=\overline{\mathcal{P}}_{\mu,\varphi}^{q,t}(F_i\cup\,H_i)
=\overline{\mathcal{P}}_{\mu,\varphi}^{q,t}(F_i)+\overline{\mathcal{P}}_{\mu,\varphi}^{q,t}(H_i).
$$
Hence,
$$
\mathcal{P}_{\mu,\varphi}^{q,t}(A)+\mathcal{P}_{\mu,\varphi}^{q,t}(B)
\leq\displaystyle\sum_i\overline{\mathcal{P}}_{\mu,\varphi}^{q,t}(E_i)\leq\mathcal{P}_{\mu,\varphi}^{q,t}(A\cup\,B)+\varepsilon
$$
and the result is obtained by having $\varepsilon\downarrow0$.
\begin{lemma}\label{borelhmuphimesurable}
	Borel sets are $\mathcal{P}_{\mu,\varphi}^{q,t}$-measurable.
\end{lemma}
\hskip-10pt The \textbf{Proof} is similar as in Lemma \ref{borelhmuphimesurable} and thus it is left to the reader.
\begin{definition}
	The restriction of $\mathcal{P}_{\mu,\varphi}^{q,t}$ on Borel sets is called the mixed generalized packing measure on $\mathbb{R}^d$.
\end{definition}
It holds as for the case of the multifractal analysis of a single measure that each of the measures $\mathcal{H}_{\mu}^{q,t}$, $\mathcal{P}_{\mu,\varphi}^{q,t}$ and the pre-measure $\overline{\mathcal{P}}_{\mu,\varphi}^{q,t}$ assign a dimension to every set $E\subseteq\mathbb{R}^d$.
\begin{proposition}\label{anouarproposition1}
	Given a subset $E\subseteq\mathbb{R}^d$,
	\begin{enumerate}
		\item There exists a unique number $\hbox{dim}_{\mu,\varphi}^q(E)\in[-\infty,+\infty]$ such that
		$$
		\mathcal{H}_{\mu,\varphi}^{q,t}(E)=\left\{\begin{array}{lll}+\infty&\hbox{for}&t<\hbox{dim}_{\mu,\varphi}^q(E)\hfill\cr\medskip
		0&\hbox{si}&t>\hbox{dim}_{\mu,\varphi}^q(E)\end{array}\right.
		$$
		\item There exists a unique number $\hbox{Dim}_{\mu,\varphi}^q(E)\in[-\infty,+\infty]$ such that
		$$
		\mathcal{P}_{\mu,\varphi}^{q,t}(E)=\left\{\begin{array}{lll}+\infty&\hbox{for}&t<\hbox{Dim}_{\mu,\varphi}^q(E)\hfill\cr\medskip
		0&\hbox{for}&t>\hbox{Dim}_{\mu,\varphi}^q(E)\end{array}\right.
		$$
		\item There exists a unique number $\Delta_\mu^q(E)\in[-\infty,+\infty]$ such that
		$$
		\overline{\mathcal{P}}_{\mu,\varphi}^{q,t}(E)=\left\{\begin{array}{lll}+\infty&\hbox{for}&t<\Delta_\mu^q(E)\hfill\cr\medskip
		0&\hbox{for}&t>\Delta_\mu^q(E)\end{array}\right.
		$$
	\end{enumerate}
\end{proposition}
\begin{definition}\label{mixeddimensions}
	The quantities $\hbox{dim}_{\mu,\varphi}^q(E)$, $\hbox{Dim}_{\mu,\varphi}^q(E)$ and $\Delta_\mu^q(E)$ define the so-called mixed multifractal generalizations of the Hausdorff dimension, the packing dimension and the logarithmic index of the set $E$.
\end{definition}
Remark that for $k=1$ and $\varphi$ the log function $\varphi(r)=\log(r)$, we come back to the classical definitions of the Hausdorff and packing measures and dimensions in their original forms (by taking $q=0$) and their generalized multifractal variants for $q$ being arbitrary. The mixed case studied here may be also applied for a single measure and thus the results and characterizations outpointed in the present work remains valid for a single measure. Indeed, denote $Q_i=(0,0,...,q_i,0,...,0)$ the vector with zero coordinates except the ith one which equals $q_i$, we obtain the multifractal generalizations of the Hausdorff $\varphi$-measure and $\varphi$-dimension, the packing $\varphi$-dimension and the logarithmic $\varphi$-index of the set $E$ for the single measure $\mu_i$,
$$
\hbox{dim}_{\mu,\varphi}^{Q_i}(E)=\hbox{dim}_{\mu_i,\varphi}^{q_i}(E),
$$
$$
\hbox{Dim}_{\mu,\varphi}^{Q_i}(E)=\hbox{Dim}_{\mu_i,\varphi}^{q_i}(E)
$$
and
$$
\Delta_{\mu,\varphi}^{Q_i}(E)=\Delta_{\mu_i,\varphi}^{q_i}(E).
$$
Similarly, for the null vector of $\mathbb{R}^k$, we obtain
$$
\hbox{dim}_{\mu,\varphi}^{0}(E)=\hbox{dim}_\varphi(E),
$$
$$
\hbox{Dim}_{\mu,\varphi}^{0}(E)=\hbox{Dim}_\varphi(E)
$$
and
$$
\Delta_{\mu,\varphi}^{0}(E)=\Delta_\varphi(E).
$$
We may obtain further
$$
\hbox{dim}_{\mu,\log}^{Q_i}(E)=\hbox{dim}_{\mu_i,\log}^{q_i}(E),
$$
$$
\hbox{Dim}_{\mu,\log}^{Q_i}(E)=\hbox{Dim}_{\mu_i,\log}^{q_i}(E)
$$
and
$$
\Delta_{\mu,\log}^{Q_i}(E)=\Delta_{\mu_i,\log}^{q_i}(E).
$$
Similarly, for the null vector of $\mathbb{R}^k$, we obtain
$$
\hbox{dim}_{\mu,\log}^{0}(E)=\hbox{dim}_{\log}(E)=\hbox{dim}(E),
$$
$$
\hbox{Dim}_{\mu,\log}^{0}(E)=\hbox{Dim}_{\log}(E)=\hbox{Dim}(E)
$$
and
$$
\Delta_{\mu,\log}^{0}(E)=\Delta_{\log}(E)=\Delta(E).
$$
\textbf{Proof of Proposition \ref{anouarproposition1}.} We will sketch only the proof of the first point. The rest is analogous.\\
1. We claim that $\forall\,t\in\mathbb{R}$ such that $\mathcal{H}_{\mu,\varphi}^{q,t}(E)<\infty$ it holds that 
$$
\mathcal{H}_{\mu,\varphi}^{q,t'}(E)=0\,,\;\;\forall\,t'>t.
$$
Indeed, let $0<\epsilon<1$, $F\subseteq\,E$ and $(B(x_i,r_i))_i$ be a centered $\epsilon$-covering of $F$. We have
$$
\begin{array}{lll}
\medskip\overline{\mathcal{H}}_{\mu,\epsilon}^{q,t'}(F)
&\leq&\displaystyle\sum_i(\mu(B(x_i,r_i)))^qe^{t'\varphi(r_i)}\hfill\cr\medskip
&\leq&e^{(t'-t)\varphi(\varepsilon)}\displaystyle\sum_i(\mu(B(x_i,r_i)))^qe^{t\varphi(r_i)}.\hfill
\end{array}
$$
Consequently,
$$
{\overline{H}}_{\mu,\varphi,\epsilon}^{q,t'}(F)\leq\,e^{(t'-t)\varphi(\varepsilon)}{\overline{H}}_{\mu,\varphi,\epsilon}^{q,t}(F).
$$
Hence, as $\varphi(\varepsilon)\rightarrow-\infty$ as $\varepsilon\rightarrow0$, we obtain
$$
\overline{\mathcal{H}}_{\mu,\varphi}^{q,t'}(F)=0,\quad\forall\,F\subseteq\,E.
$$
As a result, $\mathcal{H}_{\mu,\varphi}^{q,t'}(E)=0$. We then set
$$
\hbox{dim}_{\mu,\varphi}^q(E)=\inf\{\,t\in\mathbb{R};\;\;\mathcal{H}_{\mu,\varphi}^{q,t}(E)=0\,\}.
$$
One can proceed otherwise by claiming that $\forall\,t\in\mathbb{R}$ such that $\mathcal{H}_{\mu,\varphi}^{q,t}(E)>0$ it holds that 
$$
\mathcal{H}_{\mu,\varphi}^{q,t'}(E)=+\infty\,,\;\;\forall\,t'<t.
$$
Indeed, proceeding as previously, we obtain for $\epsilon>0$,
$$
e^{(t'-t)\varphi(\epsilon)}{\overline{H}}_{\mu,\epsilon}^{q,t}(F)\leq{\overline{H}}_{\mu,\epsilon}^{q,t'}(F).
$$
Hence, for the same reasons as above,
$$
\overline{\mathcal{H}}_{\mu,\varphi}^{q,t'}(F)=+\infty,\quad\forall\,F\subseteq\,E.
$$
As a result, $\mathcal{H}_{\mu,\varphi}^{q,t'}(E)=+\infty$. We then set
$$
\hbox{dim}_{\mu}^q(E)=\sup\{\,t\in\mathbb{R};\;\;\mathcal{H}_{\mu,\varphi}^{q,t}(E)=+\infty\,\}.
$$
Next, we aim to study the characteristics of the mixed multifractal generalizations of dimensions. We now adapt the following notations. For $q=(q_1,...,q_k)\in\mathbb{R}^k$,
$$
b_{\mu,\varphi}(q,E)=\hbox{dim}_{\mu,\varphi}^{q}(E),
$$
$$
B_{\mu,\varphi}(q,E)=\hbox{Dim}_{\mu,\varphi}^{q}(E)
$$
and
$$
\Lambda_{\mu,\varphi}(q,E)=\Delta_{\mu,\varphi}^{q}(E).
$$
When $E=\hbox{supp}(\mu)$ is the support of the measure $\mu$, we will omit the indexation with $E$ and denote simply
$$
b_{\mu,\varphi}(q),\,\,B_{\mu,\varphi}(q)\,\,\hbox{and}\,\,\Lambda_{\mu,\varphi}(q).
$$
The following propositions resume the characteristics of these functions and extends the results of L. Olsen \cite{olsen1} to our case.
\begin{proposition}\label{anouarproposition2} The following assertions hold.
	\begin{description}
		\item[a.] $b_{\mu,\varphi}(q,.)$ and $B_{\mu,\varphi}(q,.)$ are non decreasing with respect to the inclusion property in $\mathbb{R}^d$.
		\item[b.] $b_{\mu,\varphi}(q,.)$ and $B_{\mu,\varphi}(q,.)$ are $\sigma$-stable.
	\end{description}
\end{proposition}
\hskip-10pt\textbf{Proof.}
{\bf a.} Let $E\subseteq\,F$ be subsets of $\mathbb{R}^d$. We have
$$
\mathcal{H}_{\mu,\varphi}^{q,t}(E)=\displaystyle\sup_{A\subseteq\,E}\overline{\mathcal{H}}_{\mu,\varphi}^{q,t}(A)
\leq\displaystyle\sup_{A\subseteq\,F}\overline{\mathcal{H}}_{\mu,\varphi}^{q,t}(A)=\mathcal{H}_{\mu,\varphi}^{q,t}(F).
$$
So for the monotony of $b_{\mu,\varphi}(q,.)$. \\
{\bf b.} Let $(A_n)_n$ be a countable set of subsets $A_n\subseteq\mathbb{R}^d$ and denote $A=\displaystyle\bigcup_nA_n$. It holds from the monotony of $b_{\mu,\varphi}(q,.)$ that
$$
b_{\mu,\varphi}(q,A_n)\leq\,b_{\mu,\varphi}(q,A),\quad\forall\,n.
$$
Hence,
$$
\displaystyle\sup_nb_{\mu,\varphi}(q,A_n)\leq\,b_{\mu,\varphi}(q,A).
$$
Next, for any $t>\displaystyle\sup_nb_{\mu,\varphi}(q,A_n)$, there holds that
$$
\mathcal{H}_{\mu,\varphi}^{q,t}(A_n)=0,\quad\forall\,n.
$$
Consequently, from the sub-additivity property of $\mathcal{H}_{\mu,\varphi}^{q,t}$, it holds that
$$
\mathcal{H}_{\mu,\varphi}^{q,t}(\displaystyle\bigcup_nA_n)=0,\quad\forall\,t>\displaystyle\sup_nb_{\mu,\varphi}(q,A_n).
$$
Which means that
$$
b_{\mu,\varphi}(q,A)\leq\,t,\quad\forall\,t>\displaystyle\sup_nb_{\mu,\varphi}(q,A_n).
$$
Hence,
$$
b_{\mu,\varphi}(q,A)\leq\displaystyle\sup_nb_{\mu,\varphi}(q,A_n).
$$
Similar arguments permit to prove the properties of $B_{\mu,\varphi}(q,.)$.

Next, we continue to study the characteristics of the mixed generalized multifractal $\varphi$-dimensions. The following result is obtained.
\begin{proposition}\label{anouarproposition3} The following assertions are true.
	\begin{description}
		\item[a.] The functions $q\longmapsto\,B_{\mu,\varphi}(q)$ and $q\longmapsto\Lambda_{\mu,\varphi}(q)$ are convex.
		\item[b.] For $i=1,2,...,k$ and $\widehat{q_i}=(q_1,\dots,q_{i-1},q_{i+1},\dots,q_k)$ fixed, the functions $q_i\longmapsto\,b_{\mu,\varphi}(q)$, $q_i\longmapsto\,B_{\mu,\varphi}(q)$ and $q_i\longmapsto\Lambda_{\mu,\varphi}(q)$ are non increasing.
	\end{description}
\end{proposition}
\hskip-10pt{\bf Proof.}
{\bf a.} We start by proving that $\Lambda_{\mu,,\varphi}(.,E)$ is convex. Let $p,\,q\in\mathbb{R}^k$, $\alpha\in]0,1[$ and let also
$$
s>\Lambda_{\mu,\varphi}(p,E)\;\hbox{and}\;t>\Lambda_{\mu,\varphi}(q,E).
$$
Consider next a centered $\epsilon$-packing $(B_i=B(x_i,r_i))_i$ of $E$. Applying H\"older's inequality, it holds that
$$
\begin{array}{lll}
\medskip&&\displaystyle\sum_i(\mu(B_i))^{\alpha\,q+(1-\alpha)p}e^{(\alpha\,t+(1-\alpha)s)\varphi(r_i)}\hfill\cr\medskip
&\leq&\biggl(\displaystyle\sum_i(\mu(B_i))^qe^{t\varphi(r_i)}\biggr)^\alpha\biggl(\displaystyle\sum_i(\mu(B_i))^pe^{s\varphi(r_i)}\biggr)^{1-\alpha}.
\end{array}
$$
Hence,
$$
\overline{\mathcal{P}}_{\mu,\epsilon}^{\alpha\,q+(1-\alpha)p,\alpha\,t+(1-\alpha)s}(E)\leq\biggl(\overline{\mathcal{P}}_{\mu,\epsilon}^{q,t}(E)\biggr)^\alpha
\biggl(\overline{\mathcal{P}}_{\mu,\epsilon}^{p,s}(E)\biggr)^{1-\alpha}.
$$
The limit on $\epsilon\downarrow0$ gives
$$
\overline{\mathcal{P}}_{\mu,\varphi}^{\alpha q+(1-\alpha)p,\alpha\,t+(1-\alpha)s}(E)\leq\biggl(\overline{\mathcal{P}}_{\mu,\varphi}^{q,t}(E)\biggr)^\alpha
\biggl(\overline{\mathcal{P}}_{\mu,\varphi}^{p,s}(E)\biggr)^{1-\alpha}.
$$
Consequently,
$$
\overline{\mathcal{P}}_{\mu,\varphi}^{\alpha q+(1-\alpha)p,\alpha\,t+(1-\alpha)s}(E)=0\,,\;\;\forall\,s>\Lambda_{\mu,E}(p)\;\;\hbox{and}\;\;t>\Lambda_{\mu,E}(q).
$$
It results that
$$
\Lambda_{\mu\varphi}(\alpha q+(1-\alpha)p,E)\leq\alpha\Lambda_{\mu,\varphi}(q,E)+(1-\alpha)\Lambda_{\mu,\varphi}(p,E).
$$
We now prove the convexity of $B_{\mu,\varphi}(.,E)$. We set in this case
$$
t=B_{\mu,\varphi}(q,E)\;\;\hbox{and}\;\;s=B_{\mu,\varphi}(p,E).
$$
We have
$$
\mathcal{P}_{\mu,\varphi}^{q,t+\varepsilon}(E)=\mathcal{P}_{\mu,\varphi}^{p,s+\varepsilon}(E)=0.
$$
Therefore, there exists $(H_i)_i$ and $(K_i)_i$ coverings of the set $E$ for which
$$
\displaystyle\sum_i\overline{\mathcal{P}}_{\mu,\varphi}^{q,t+\varepsilon}(H_i)\leq1
\qquad\hbox{et}\qquad
\displaystyle\sum_i\overline{\mathcal{P}}_{\mu,\varphi}^{p,s+\varepsilon}(K_i)\leq1.
$$
Denote for $n\in\mathbb{N}$, $E_n=\displaystyle\bigcup_{1\leq\,i,j\leq\,n}(H_i\cap K_j)$. Thus, $(E_n)_n$ is a covering of $E$. So that,
$$
\begin{array}{lll}& &\mathcal{P}_{\mu,\varphi}^{\alpha\,q+(1-\alpha)p,\alpha\,t+(1-\alpha)s+\varepsilon}(E_n)\hfill\cr\medskip
&\leq&\displaystyle\sum_{i,j=1}^n\mathcal{P}_{\mu,\varphi}^{\alpha\,q+(1-\alpha)p,\alpha t+(1-\alpha)s+\varepsilon}(H_i\cap\,K_j)\hfill\cr\medskip
&\leq&\displaystyle\sum_{i,j=1}^n\overline{\mathcal{P}}_{\mu,\varphi}^{\alpha\,q+(1-\alpha)p,\alpha t+(1-\alpha)s+\varepsilon}(H_i\cap\,K_j)\hfill\cr\medskip
&\leq&\biggl(\displaystyle\sum_{i,j=1}^n\overline{\mathcal{P}}_{\mu,\varphi}^{q,t+\varepsilon}(H_i\cap\,K_j)\biggr)^\alpha
\biggl(\displaystyle\sum_{i,j=1}^n\overline{\mathcal{P}}_{\mu,\varphi}^{p,s+\varepsilon}(H_i\cap\,K_j)\biggr)^{1-\alpha}\hfill\cr\medskip
&\leq&n^\alpha\,n^{1-\alpha}=n<\infty.\end{array}
$$
Consequently,
$$
B_{\mu,\varphi}(\alpha\,q+(1-\alpha)p,E_n)\leq\alpha\,t+(1-\alpha)s+\varepsilon,\quad\forall\,\varepsilon>0.
$$
Hence,
$$
B_{\mu,\varphi}(\alpha\,q+(1-\alpha)p,E)\leq\alpha\,B_{\mu,\varphi}(q,E)+(1-\alpha)B_{\mu,\varphi}(p,E).
$$
{\bf b.} For $i=1,2,\dots,k$, let $\widehat{q_i}$ fixed and $p_i\leq\,q_i$ reel numbers. Denote next
$$
q=(q_1,\dots,q_{i-1},q_i,q_{i+1},\dots,q_k)\;\;\hbox{and}\;\;p=(q_1,\dots,q_{i-1},p_i,q_{i+1},\dots,q_k).
$$
Let finally $A\subseteq\,E$. For a centered $\epsilon$-covering $(B(x_i,r_i))_i$ of $A$, we have immediately
$$
\mu(B(x_i,r_i))^{q}e^{t\varphi(r_i)}\leq\mu(B(x_i,r_i))^{p}e^{t\varphi(r_i)},\;\;\forall\,t\in\mathbb{R}.
$$
Hence,
$$
{\overline{H}}_{\mu,\varphi,\epsilon}^{q,t}(A)\leq{\overline{H}}_{\mu,\varphi,\epsilon}^{p,t}(A).
$$
When $\epsilon\downarrow0$, we obtain
$$
{\overline{H}}_{\mu,\varphi}^{q,t}(A)\leq{\overline{H}}_{\mu,\varphi}^{p,t}(A).
$$
Therefore,
$$
\mathcal{H}_{\mu,\varphi}^{q,t}(E)=\displaystyle\sup_{A\subseteq\,E}\overline{\mathcal{H}}_{\mu,\varphi}^{q,t}(A)
\leq\displaystyle\sup_{A\subseteq\,E}\overline{\mathcal{H}}_{\mu,\varphi}^{p,t}(A)={\mathcal{H}}_{\mu,\varphi}^{p,t}(E).
$$
This induces the fact that
$$
{\mathcal{H}}_{\mu,\varphi}^{q,t}(E)=0,\quad\forall\,t>b_{\mu,\varphi}(p,E).
$$
Consequently
$$
b_{\mu,\varphi}(q,E)<t,\quad\forall\,t>b_{\mu,\varphi}(p,E).
$$
Hence,
$$
b_{\mu,\varphi}(q,E)\leq\,b_{\mu,\varphi}(p,E).
$$
The remaining part to prove the monotony $\Lambda_{\mu,\varphi}(.,E)$ and $B_{\mu,\varphi}(.,E)$ is analogous.
\begin{proposition}\label{anouarproposition4} The following assertions are true.
	\begin{description}
		\item[a.] $0\leq\,b_{\mu,\varphi}(q)\leq\,B_{\mu,\varphi}(q)\leq\Lambda_{\mu,\varphi}(q)$, whenever $q_i<1$ for all $i=1,2,...,k$.
		\item[b.] $b_{\mu,\varphi}(\P_i)=B_{\mu,\varphi}(\P_i)=\Lambda_{\mu,\varphi}(\P_i)=0$, where $\P_i=(0,0,...,1,0,...,0)$.
		\item[c.] $b_{\mu,\varphi}(q)\leq\,B_{\mu,\varphi}(q)\leq\Lambda_{\mu,\varphi}(q)\leq0$ whenever $q_i>1$ for all $i=1,2,...,k$.
	\end{description}
\end{proposition}
\hskip-10pt The proof of this results reposes on the following intermediate ones.
\begin{lemma}\label{anouarlemme3}
	There exists a constant $\xi\in]0,+\infty[$ satisfying for any $E\subseteq\mathbb{R}^d$,
	$$
	\mathcal{H}_{\mu,\varphi}^{q,t}(E)\leq\xi\mathcal{P}_{\mu,\varphi}^{q,t}(E)\leq\xi\overline{\mathcal{P}}_{\mu,\varphi}^{q,t}(E),\qquad\forall\,q,t.
	$$
	More precisely, $\xi$ is the number related to the Besicovitch covering theorem.
\end{lemma}
\hskip-10pt{\bf Proof.} It suffices to prove the first inequality. The second is always true for all $\xi>0$. Let $F\subseteq\mathbb{R}^d$, $\epsilon>0$ and $\mathcal{V}=\{\,B(x,\frac{\epsilon}{2});\;\;\,x\in\,F\,\}$. Let next $\bigl((B_{ij})_j\bigr)_{1\leq\,i\leq\xi}$ be the $\xi$ sets of $\mathcal{V}$ obtained by the Besicovitch covering theorem. So that, $(B_{ij})_{i,j}$ is a centered $\epsilon$-covering of the set $F$ and for each $i$, $(B_{ij})_j$ is a centered $\epsilon$-packing of $F$. Therefore,
$$
\displaystyle\overline{\mathcal{H}}_{\mu,\varphi,\epsilon}^{q,t}(F)\leq\sum_{i=1}^{\xi}\sum_j\bigl(\mu(B_{ij})\bigr)^qe^{t\varphi(r_{ij})}
\leq\sum_{i=1}^{\xi}\overline{\mathcal{P}}_{\mu,\varphi,\epsilon}^{q,t}(F)=\xi\overline{\mathcal{P}}_{\mu,\varphi,\epsilon}^{q,t}(F).
$$
Hence,
$$
\overline{\mathcal{H}}_{\mu,\varphi}^{q,t}(F)\leq\xi\overline{\mathcal{P}}_{\mu,\varphi}^{q,t}(F).
$$
Consequently, for $E\subseteq\displaystyle\bigcup_iE_i$, we obtain
$$
\begin{array}{lll}
\mathcal{H}_{\mu,\varphi}^{q,t}(E)&=&\mathcal{H}_{\mu,\varphi}^{q,t}(\displaystyle\bigcup_i(E_i\cap\,E))\hfill\cr\medskip
&\leq&\displaystyle\sum_i\mathcal{H}_{\mu,\varphi}^{q,t}(E_i\cap\,E)\hfill\cr\medskip
&\leq&\displaystyle\sum_i\displaystyle\sup_{F\subseteq E_i\cap\,E}\overline{\mathcal{H}}_{\mu,\varphi}^{q,t}(F)\hfill\cr\medskip
&\leq&\xi\displaystyle\sum_i\displaystyle\sup_{F\subseteq E_i\cap\,E}\overline{\mathcal{P}}_{\mu,\varphi}^{q,t}(F)\hfill\cr\medskip
&\leq&\xi\displaystyle\sum_i\overline{\mathcal{P}}_{\mu,\varphi}^{q,t}(E_i).\end{array}
$$
So as Lemma \ref{anouarlemme3}.\\
{\bf Proof of Proposition \ref{anouarproposition4}.} It follows from Proposition \ref{anouarproposition2}, Proposition \ref{anouarproposition3} and Lemma \ref{anouarlemme3}.
\section{$\varphi$-mixed multifractal generalization of Bouligand-Minkowski's dimensions}
Let $k\geq1$ be en integer and  $\mu=(\mu_1,\,\mu_2,\dots,\mu_k)$ be a vector valued measure composed of probability measures on $\mathbb{R}^d$. Denote as previously
$$
\mu(B(x,r))\equiv\bigl(\mu_1(B(x,r)),\dots,\,\mu_k(B(x,r))\bigr)
$$
and for $q=(q_1,q_2,\dots,q_k)\in\mathbb{R}^k$,
$$
(\mu(B(x,r)))^q\equiv(\mu_1(B(x,r)))^{q_1}\dots(\mu_k(B(x,r)))^{q_k}.
$$
Next, for a nonempty subset $E\subseteq\hbox{supp}(\mu)$, $\delta>0$ and $q\in\mathbb{R}$, we put
$$
\mathcal{T}_{\mu,\delta}^q(E)=\inf\left\{\displaystyle\sum_i\bigl(\mu\bigl(B(x_i,\delta)\bigr)\bigr)^q\right\}
$$
where the inf is over the set of all centered $\delta$-coverings $\bigl(B(x_i,\delta)\bigr)_i$ of the set $E$. The mixed multifractal generalized Bouligand-Minkowski $\varphi$-dimensions are
$$
{\overline{L}}_{\mu,\varphi}^q(E)=\displaystyle\limsup_{\delta\downarrow0}
\displaystyle\frac{\log\bigl(\mathcal{T}_{\mu,\delta}^q(E)\bigr)}{-\varphi(\delta)}
$$
for the upper one and
$$
{\underline{L}}_{\mu,\varphi}^q(E)=\displaystyle\liminf_{\delta\downarrow0}
\displaystyle\frac{\log\bigl(\mathcal{T}_{\mu,\delta}^q(E)\bigr)}{-\varphi(\delta)}
$$
for the lower. In the case of equality, the common value is denoted ${L}_{\mu,\varphi}^q(E)$ and is called the mixed multifractal generalized Bouligand-Minkowski $\varphi$-dimension of the set $E$.

Such dimensions may also be obtained via the $\delta$-packings as follows. Indeed, for $\delta>0$ and $q\in\mathbb{R}$, we set
$$
\mathcal{S}_{\mu,\delta}^q(E)=\sup\left\{\displaystyle\sum_i\bigl(\mu\bigl(B(x_i,\delta)\bigr)\bigr)^q\right\}
$$
where the sup is taken over all the centered $\delta$-packings $\bigl(B(x_i,\delta)\bigr)_i$ of the set $E$. The upper dimension is
$$
{\overline{C}}_{\mu,\varphi}^q(E)=\displaystyle\limsup_{\delta\downarrow0}
\displaystyle\frac{\log\bigl(\mathcal{S}_{\mu,\delta}^q(E)\bigr)}{-\varphi(\delta)}
$$
and the lower is
$$
{\underline{C}}_{\mu,\varphi}^q(E)=\displaystyle\liminf_{\delta\downarrow0}
\displaystyle\frac{\log\bigl(\mathcal{S}_{\mu,\delta}^q(E)\bigr)}{-\varphi(\delta)}
$$
and similarly, when these are equal, the common value will be denoted ${C}_{\mu,\varphi}^q(E)$ and it defines the dimension of $E$.
\begin{definition}\label{bouligandminkowskimixte}
	For $E\subseteq\hbox{supp}(\mu)$ and $q=(q_1,q_2,\dots,q_k)\in\mathbb{R}^k$, we will call
	\begin{description}
		\item[a.] ${\overline{C}}_{\mu,\varphi}^q(E)$ and ${\overline{L}}_{\mu,\varphi}^q(E)$ the upper $\mu$-mixed multifractal generalized Bouligand Minkowski $\varphi$-dimension of $E$.
		\item[b.] ${\underline{C}}_{\mu,\varphi}^q(E)$ and ${\underline{L}}_{\mu,\varphi}^q(E)$ the lower $\mu$-mixed multifractal generalized Bouligand Minkowski $\varphi$-dimension of $E$.
		\item[c.] ${{C}}_{\mu,\varphi}^q(E)$ and ${{L}}_{\mu,\varphi}^q(E)$ the $\mu$-mixed multifractal generalized Bouligand Minkowski $\varphi$-dimension of $E$.
	\end{description}
\end{definition}
\begin{remark} As for the classical multifractal contexts, each of the quantities above defines in fact a mixed generalization that can be different from the other.
\end{remark}
\begin{theorem}\label{anouartheorem1}
	\begin{enumerate}
		\item For all $q\in\mathbb{R}^{k}$, we have
		$$
		{\underline{L}}_{\mu,\varphi}^q(E)\leq{\underline{C}}_{\mu,\varphi}^q(E)\quad\hbox{and}\quad
		{\overline{L}}_{\mu,\varphi}^q(E)\leq{\overline{C}}_{\mu,\varphi}^q(E).
		$$
		\item For any $q\in\mathbb{R}^{*\,k}_-$, we have
		\begin{description}
			\item[i.] $b_{\mu,\varphi}(q,E)\leq{\underline{L}}_{\mu,\varphi}^q(E)={\underline{C}}_{\mu,\varphi}^q(E)$.
			\item[ii.] ${\overline{L}}_{\mu,\varphi}(q,E)={\overline{C}}_{\mu,\varphi}^q(E)={\Lambda}_{\mu,\varphi}(q,E)$.
		\end{description}
		\item For any $q\in\mathbb{R}^{*\,k}_+$, we have
		$$
		{\overline{L}}_{\mu,\varphi}(q,E)\leq{\overline{C}}_{\mu,\varphi}^q(E)\leq{\Lambda}_{\mu,\varphi}(q,E).
		$$
	\end{enumerate}
\end{theorem}
\hskip-10pt{\bf Proof.} {\bf 1.} Using Besicovitch covering theorem we get
$$
\mathcal{T}_{\mu,\varphi,\delta}^q(E)\leq\,C\mathcal{S}_{\mu,\varphi,\delta}^q(E),
$$
with some constant $C$ fixed. So {\bf 1.} is proved.\\
{\bf 2.} We firstly prove that
$$
{\underline{L}}_{\mu,\varphi}^q(E)\geq{\underline{C}}_{\mu,\varphi}^q(E)\quad\hbox{and}\quad
{\overline{L}}_{\mu,\varphi}^q(E)\geq{\overline{C}}_{\mu,\varphi}^q(E).
$$
Indeed, let $\bigl(B(x_i,\delta)\bigr)_i$ be a centered $\delta$-packing of $E$ and $\bigl(B(y_i,\frac{\delta}{2})\bigr)$ be a centered $\frac{\delta}{2}$-covering of $E$. Consider for each $i$, the integer $k_i$ such that $x_i\in B(y_{k_i},\frac{\delta}{2})$. It is straightforward that for $i\not=j$ we have $k_i\not=k_j$. Consequently, for $q\in\mathbb{R}^{*\,k}_-$, there holds that
$$
\begin{array}{lll}\displaystyle\sum_i\bigl(\mu(B(x_i,\delta))\bigr)^q
&=&\displaystyle\sum_i\Biggl(\displaystyle\frac{\mu(B(x_i,\delta))}{\mu(B(y_{k_i},\delta/2))}\Biggr)^q
\bigl(\mu(B(y_{k_i},\frac{\delta}{2}))\bigr)^q\hfill\cr\medskip
&\leq&\displaystyle\sum_i\bigl(\mu(B(y_i,\frac{\delta}{2}))\bigr)^q.\end{array}
$$
Which means that
$$
\mathcal{S}_{\mu,\varphi,\delta}^q(E)\leq\mathcal{T}_{\mu,\varphi,\frac{\delta}{2}}^q(E)
$$
and thus, for any $q\in\mathbb{R}^{*\,k}_-$,
$$
{\underline{L}}_{\mu,\varphi}^q(E)\geq{\underline{C}}_{\mu,\varphi}^q(E)\quad\hbox{and}\quad
{\overline{L}}_{\mu,\varphi}^q(E)\geq{\overline{C}}_{\mu,\varphi}^q(E)
$$
Using assertion {\bf 1.}, we obtain the equalities
$$
{\underline{L}}_{\mu,\varphi}^q(E)={\underline{C}}_{\mu,\varphi}^q(E)\quad\hbox{and}\quad
{\overline{L}}_{\mu,\varphi}^q(E)={\overline{C}}_{\mu,\varphi}^q(E)
$$
for all $q\in\mathbb{R}^{*\,k}_-$. Therefore, to prove {\bf 2.i.}, it remains to prove the inequality of the left hand side. So, let $t>{\underline{L}}_{\mu,\varphi}^q(E)$ and $F\subseteq E$. Consider next a sequence $(\delta_n)_n\subseteq]0,1[$ to be $\downarrow0$, and satisfying
$$
t>\displaystyle\frac{\log(\mathcal{T}_{\mu,\varphi,\delta_n}^q(E))}{-\varphi(\delta_n)},\quad\forall\,n\in\mathbb{N}.
$$
This means that for each $n\in\mathbb{N}$, there exists a centered $\delta_n$-covering $\bigl(B(x_{ni},\delta_n)\bigr)_i$ of $E$ such that
$$
\displaystyle\sum_i\bigl(\mu(B(x_{ni},\delta_n))\bigr)^q<e^{-t\varphi(\delta_n)}.
$$
There balls may be considered to be intersecting the set $F$. Next, for each $i$, choose an element $y_i\in\,B(x_{ni},\delta_n)\cap\,F$.
This results on a centered $2\delta_n$-covering $\bigl(B(y_i,2\delta_n)\bigr)_i$ of $F$. Therefore,
$$
\begin{array}{lll}\overline{\mathcal{H}}_{\mu,\varphi,2\delta_n}^{q,t}(F)
&\leq&\displaystyle\sum_i\bigl(\mu(B(x_{ni},\delta_n))\bigr)^qe^{t\varphi(2\delta_n)}\hfill\cr\medskip
&=&C_t\displaystyle\sum_i\Biggl(\displaystyle\frac{\mu(B(y_i,2\delta_n))}{\mu(B(x_{ni},\delta_n))}\Biggr)^q
\bigl(\mu(B(x_{ni},\delta_n))\bigr)^qe^{t\varphi(\delta_n)}\hfill\cr\medskip
&\leq&C_t\displaystyle\sum_i\bigl(\mu(B(x_{ni},\delta_n))\bigr)^qe^{t\varphi(\delta_n)}\hfill\cr\medskip
&\leq&C_tC'_t.\end{array}
$$
Hence,
$$
\overline{\mathcal{H}}_{\mu,\varphi}^{q,t}(F)\leq C\quad\forall\,F\subseteq\,E,\;\;t>{\underline{L}}_{\mu,\varphi}^q(E).
$$
So that,
$$
\mathcal{H}_{\mu,\varphi}^{q,t}(E)\leq C<\infty,\quad\forall\,t>{\underline{L}}_{\mu,\varphi}^q(E).
$$
Consequently,
$$
b_{\mu,\varphi}(q,E)\leq t,\quad\forall\,t>{\underline{L}}_{\mu,\varphi}^q(E)\,\Rightarrow\,b_{\mu,\varphi}(q,E)\leq{\underline{L}}_{\mu,\varphi}^q(E).
$$
The remaining part can be proved by following similar techniques.\\

Next we need to introduce the following quantities which will be useful later. Let $\mu=(\mu_1,\,\mu_2,\dots,\mu_k)$ be a vector valued measure composed of probability measures on $\mathbb{R}^d$. For $j=1,\,2,\,\dots,\,k$, $a>1$ and $E\subseteq\hbox{supp}(\mu)$, denote
$$
T_a^j(E)=\displaystyle\limsup_{r\downarrow0}\Bigl(\displaystyle\sup_{x\in\,E}\frac{\mu_j\bigl(B(x,ar)\bigr)}{\mu_j\bigl(B(x,r)\bigr)}\Bigr)
$$
and for $x\in\hbox{supp}(\mu)$, $T_a^j(x)=T_a^j(\{x\})$. Denote also
$$
P_0(\mathbb{R}^d,E)=\{\,\mu\,;\;\;\exists\,a\,>\,1\,;\;\;\forall\,x\in\,E,\;\;T_a^j(x)<\infty,\;\;\forall\,j\,\},
$$
$$
P_1(\mathbb{R}^d,E)=\{\,\mu\,;\;\;\exists\,a\,>\,1\,;\;T_a^j(E)<\infty,\;\;\forall\,j\,\},
$$
$$
P_0(\mathbb{R}^d)=P_0(\mathbb{R}^d,\hbox{supp}(\mu))\qquad\hbox{and}\qquad\,P_1(\mathbb{R}^d)=P_1(\mathbb{R}^d,\hbox{supp}(\mu)).
$$
\begin{theorem}\label{anouartheorem2}
	\begin{enumerate}
		\item For $\mu\in\,P_0(\mathbb{R}^d)$ and $q\in\mathbb{R}^{*\,k}_+$, there holds that
		$$
		b_{\mu,\varphi}(q,E)\leq{\overline{L}}_{\mu,\varphi}^q(E).
		$$
		\item For $\mu\in\,P_1(\mathbb{R}^d)$ and $q\in\mathbb{R}^{*\,k}_+$, there holds that
		\begin{description}
			\item[i.] ${\underline{L}}_{\mu,\varphi}^q(E)={\underline{C}}_{\mu,\varphi}^q(E)$.
			\item[ii.] ${\overline{L}}_{\mu,\varphi}(q,E)={\overline{C}}_{\mu,\varphi}^q(E)={\Lambda}_{\mu,\varphi}(q,E)$.
		\end{description}
	\end{enumerate}
\end{theorem}
\hskip-10pt{\bf Proof.} {\bf 1.} The vector valued measure $\mu\in\,P_0(\mathbb{R}^d)$ yields that
$$
E=\displaystyle\bigcup_{m\in\mathbb{N}}E_m
$$
where
$$
E_m=\{\,x\in\,E\,;\;\;\displaystyle\frac{\mu_j(B(x_i,4r))}{\mu_j(B(x_i,r))}<m\,,\;0<r<\frac{1}{m},\,\;\;\forall\,j\,\}.
$$
Next, remark that for $t>{\overline{L}}_{\mu,\varphi}^q(E)$ and $F\subseteq\,E_m$, there exists a sequence $(\delta_n)_n\downarrow0$ for which
$$
t<\displaystyle\frac{\log(\mathcal{T}_{\mu,\delta_n}^q(F))}{-\varphi(\delta_n)},\qquad\forall\,n\in\mathbb{N}.
$$
Therefore, there exists a centered $\delta_n$-covering $(B(x_{ni},\delta_n))_i$ of $F$ satisfying
$$
\displaystyle\sum_i\bigl(\mu(B(x_{ni},\delta_n))\bigr)^q<e^{-t\varphi(\delta_n)}.
$$
Let next $y_{ni}\in\,B(x_{ni},\delta_n)$. Then, $(B(x_{ni},2\delta_n))_i$ is a centered $2\delta_n$-covering of $F$. Hence,
$$
\begin{array}{lll}\overline{\mathcal{H}}_{\mu,\varphi,2\delta_n}^{q,t}(F)
&\leq&\displaystyle\sum_i\bigl(\mu(B(y_{ni},2\delta_n))\bigr)^qe^{t\varphi(2\delta_n)}\hfill\cr\medskip
&\leq&C_t\displaystyle\sum_i\Biggl(\displaystyle\frac{\mu(B(y_{ni},2\delta_n))}{\mu(B(x_{ni},\delta_n))}\Biggr)^q
\bigl(\mu(B(x_{ni},\delta_n))\bigr)^qe^{t\varphi(\delta_n)}\hfill\cr\medskip
&\leq&C_t\displaystyle\sum_i\Biggl(\displaystyle\frac{\mu(B(x_{ni},4\delta_n))}{\mu(B(x_{ni},\delta_n))}\Biggr)^q
\bigl(\mu(B(x_{ni},\delta_n))\bigr)^qe^{t\varphi(\delta_n)}\hfill\cr\medskip
&\leq&C_tm^{|q|}\displaystyle\sum_i\bigl(\mu(B(x_{ni},\delta_n))\bigr)^qe^{t\varphi(\delta_n)}\hfill\cr\medskip
&\leq&C_tm^{|q|}C'_t\end{array}
$$
where $|q|=q_1+q_2+\dots+q_k$. Thus,
$$
\overline{\mathcal{H}}_{\mu,\varphi}^{q,t}(F)\leq C_tC'_tm^{|q|},\quad\forall\,m,\;\;\hbox{and}\;\,F\subseteq\,E_m.
$$
Which means that
$$
\mathcal{H}_{\mu,\varphi}^{q,t}(E_m)\leq C_tC'_tm^{|q|}<\infty,\quad\forall\,m,\;\;\hbox{and}\;\,t>{\underline{L}}_{\mu,\varphi}^q(E).
$$
Consequently,
$$
b_{\mu,\varphi}(q,E_m)\leq\,t,\quad\forall\,m\;\;\hbox{and}\;\,t>{\underline{L}}_{\mu,\varphi}^q(E).
$$
Using the $\sigma$-stability of $b_{\mu,\varphi}(q,.)$ (See Proposition \ref{anouarproposition2}. c.), it results that
$$
b_{\mu,\varphi}(q,E)\leq\,t,\quad\forall\,t>{\underline{L}}_{\mu,\varphi}^q(E).
$$
As a result,
$$
b_{\mu,\varphi}(q,E)\leq{\underline{L}}_{\mu,\varphi}^q(E).
$$
Assertion {\bf 2.} is left to the reader.\\

We now re-introduce the mixed multifractal generalization of the $L^q$-dimensions called also Renyi dimensions based on integral representations. See \cite{olsen3} for more details and other results. For $q\in\mathbb{R}^{*,k}$, $\mu=(\mu_1,\,\mu_2,\,\dots,\,\mu_k)$ and $\delta>0$, we set
$$
I_{\mu,\delta}^q=\displaystyle\int_{S_{\mu}}\Bigl(\mu(B(t,\delta))\Bigr)^qd\mu(t),
$$
where, in this case,
$$
S_{\mu}=\hbox{supp}(\mu_1)\times\hbox{supp}(\mu_2)\times\,\dots\,\times\hbox{supp}(\mu_k),
$$
$$
\Bigl(\mu(B(t,\delta))\Bigr)^q=\Bigl(\mu_1(B(t_1,\delta))\Bigr)^{q_1}\,\Bigl(\mu_2(B(t_2,\delta))\Bigr)^{q_2}\dots\Bigl(\mu_k(B(t_k,\delta))\Bigr)^{q_k}
$$
and
$$
d\mu(t)=d\mu_1(t_1)\,d\mu_2(t_2)\,\dots\,d\mu_k(t_k).
$$
The $(\mu,\varphi)$-mixed multifractal generalizations of the Renyi dimensions are
$$
{\overline{I}}_{\mu,\varphi}^q=\displaystyle\limsup_{\delta\downarrow0}\displaystyle\frac{\log\,I_{\mu,\delta}^q}{-\varphi(\delta)},\quad\hbox{and}\quad
{\underline{I}}_{\mu,\varphi}^q=\displaystyle\liminf_{\delta\downarrow0}\displaystyle\frac{\log\,I_{\mu,\delta}^q}{-\varphi(\delta)}.
$$
We now propose to relate these dimensions to the quantities ${\underline{C}}_{\mu,\varphi}^{q}$, ${\overline{C}}_{\mu,\varphi}^{q}$, ${\underline{L}}_{\mu,\varphi}^{q}$, ${\overline{L}}_{\mu,\varphi}^{q}$ introduced previously.
\begin{proposition}\label{anouarproposition5} The following results hold.
	\begin{description}
		\item[a.] $\forall\,q\in\mathbb{R}^{*,\,k}_-$,
		$$
		\underline{C}_{\mu,\varphi}^{q+\mathbb{I}}(\hbox{supp}(\mu))\geq\underline{I}_{\mu,\varphi}^q\quad\,and\,
		\quad\overline{C}_{\mu,\varphi}^{q+\mathbb{I}}(\hbox{supp}(\mu))\geq\overline{I}_{\mu,\varphi}^q.
		$$
		\item[b.] $\forall\,q\in\mathbb{R}^{*,\,k}_+$,
		$$
		{\underline{C}}_{\mu,\varphi}^{q+\mathbb{I}}(supp(\mu))\leq{\underline{I}_{\mu,\varphi}^q}\qquad\,and\qquad\,
		{\overline{C}}_{\mu,\varphi}^{q+\mathbb{I}}(supp(\mu))\leq{\overline{I}_{\mu,\varphi}^q}.
		$$
		\item[c.] $\forall\,q\in\mathbb{R}^{*,\,k}$, $\mu\in\,P_1(\mathbb{R}^d)$,
		$$
		{\underline{C}}_{\mu,\varphi}^{q+\mathbb{I}}(supp(\mu))=\underline{I}_{\mu,\varphi}^q\qquad\,and\qquad\,
		{\overline{C}}_{\mu,\varphi}^{q+\mathbb{I}}(supp(\mu))=\overline{I}_{\mu,\varphi}^q.
		$$
		\item[d.] $\forall\,q\in\mathbb{R}^{*,\,k}_-$,
		$$
		\underline{I}_{\mu,\varphi}^q\leq\underline{L}_{\mu,\varphi}^{q+\mathbb{I}}(supp(\mu))\qquad\,and\qquad\,
		\overline{I}_{\mu,\varphi}^q\leq\overline{L}_{\mu,\varphi}^{q+\mathbb{I}}(supp(\mu)).
		$$
	\end{description}
\end{proposition}
\hskip-10pt{\bf Proof.} We only prove {\bf a.} The remaining proofs of points {\bf b.}, {\bf c.} and {\bf d.} follow the same ideas. For $\delta>0$, let $\Bigl(B(x_i,\delta)\Bigr)_i$ be a centered $\delta$-covering of $\hbox{supp}(\mu)$ and let next $\Bigl(B(x_{ij},\delta)\Bigr)_j$, ${1\leq\,i\leq\xi}$ the $\xi$ sets defined in Besicovitch covering theorem. It holds that
$$
\begin{array}{lll}\displaystyle\sum_{i,j}\Bigl(\mu(B(x_{ij},\delta))\Bigr)^{q+\mathbb{I}}
&=&\displaystyle\sum_{i,j}\Bigl(\mu(B(x_{ij},\delta))\Bigr)^q\displaystyle\int_{B(x_{ij},\delta)^k}d\mu(t)\hfill\cr\medskip
&\geq&\displaystyle\sum_{i,j}\displaystyle\int_{B(x_{ij},\delta)^k}\Bigl(\mu(B(t,2\delta))\Bigr)^qd\mu(t)\hfill\cr\medskip
&\geq&\displaystyle\int_{S_{\mu}}\Bigl(\mu(B(t,2\delta))\Bigr)^qd\mu(t).\end{array}
$$
As a results,
$$
\xi\mathcal{S}_{\mu,\delta}^{q+\mathbb{I}}(\hbox{supp}(\mu))\geq\,I_{\mu,2\delta}^q.
$$
Which implies that
$$
\underline{C}_{\mu,\varphi}^{q+\mathbb{I}}(\hbox{supp}(\mu))\geq\underline{I}_{\mu}^q\quad\hbox{and}
\quad\overline{C}_{\mu,\varphi}^{q+\mathbb{I}}(\hbox{supp}(\mu))\geq\overline{I}_{\mu}^q.
$$
\section{A $\varphi$-mixed multifractal formalism for vector valued measures}
Let $\mu=(\mu_1,\,\mu_2,\,\dots,\,\mu_k)$ be a vector valued probability measure on $\mathbb{R}^d$. For $x\in\mathbb{R}^d$ and $j=1,2,\dots,k$, we denote
$$
{\underline\alpha}_{\mu_j}^{\varphi}(x)=\displaystyle\liminf_{r\downarrow0}\displaystyle\frac{\log(\mu_j(B(x,r)))}{\varphi(r)}\;\;\hbox{and}\;\;
{\overline\alpha}_{\mu_j}^{\varphi}(x)=\displaystyle\limsup_{r\downarrow0}\displaystyle\frac{\log(\mu_j(B(x,r)))}{\varphi(r)}
$$
respectively the local lower $\varphi$-dimension and the local upper $\varphi$-dimension of $\mu_j$ at the point $x$ and as usually the local dimension $\alpha_{\mu_j}^{\varphi}(x)$ of $\mu_j$ at $x$ will be the common value when these are equal. Next for $\alpha=(\alpha_1,\,\alpha_2,\,\dots,\,\alpha_k)\in\mathbb{R}_+^k$, let
$$
\underline{X}_{\alpha}({\varphi})=\{\,x\in\hbox{supp}(\mu)\,;\,\,{\underline\alpha}_{\mu_j}^{\varphi}(x)\geq\alpha_j\,,\forall\,j=1,2,\dots,k\,\},
$$
$$
\overline{X}^{\alpha}({\varphi})=\{\,x\in\hbox{supp}(\mu)\,;\,\,{\overline\alpha}_{\mu_j}^{\varphi}(x)\leq\alpha_j\,,\forall\,j=1,2,\dots,k\,\}
$$
and
$$
X(\alpha,{\varphi})=\underline{X}_{\alpha}({\varphi})\cap\overline{X}^{\alpha}({\varphi}).
$$
The $(\mu,{\varphi})$-mixed multifractal spectrum of the vector valued measure $\mu$ is defined by
$$
\alpha\,\longmapsto\,\hbox{dim}\,X(\alpha,\varphi)
$$
where $\hbox{dim}$ stands for the Hausdorff dimension.

In this section, we propose to compute such a spectrum for some cases of measures that resemble to the situation raised by Olsen in \cite{olsen1} but in the mixed case. This will permit to describe better the simultaneous behavior of finitely many measures. We intend precisely to compute the mixed spectrum based on the mixed multifractal generalizations of the Haudorff and packing dimensions $b_\mu$, $B_\mu$ and $\Lambda_\mu$. We start with the following technic results.
\begin{lemma}\label{spectrelemme1}
	Let $\varphi:\mathbb{R}_+\rightarrow\mathbb{R}$ be such that
	\begin{equation}\label{condition-sur-phi}
		\varphi\;\hbox{is non-decreasing and},\;\;\varphi(r)=o(\log r)\;\hbox{as}\;r\rightarrow0.
	\end{equation}
	The following assertions hold.
	\begin{description}
		\item[1.] $\forall\,\delta>0,\,t\in\mathbb{R}$ and $q\in\mathbb{R}^k_+$, $\alpha\in\mathbb{R}^k$ such that $\langle\alpha,q\rangle+t\geq0$, we have
		\begin{description}
			\item[i.] $\mathcal{H}^{\langle\alpha,q\rangle+t+k\delta}(\overline{X}^\alpha)\leq\,C\mathcal{H}_{\mu,\varphi}^{q,t}(\overline{X}^\alpha).$
			\item[ii.]
			$\mathcal{P}^{\langle\alpha,q\rangle+t+k\delta}({\overline{X}}^\alpha)\leq\,C\mathcal{P}_{\mu,\varphi}^{q,t}({\overline{X}}^\alpha).$
		\end{description}
		\item[2.] $\forall\,\delta>0,\,t\in\mathbb{R}$ and $q\in\mathbb{R}^k_-$, $\alpha\in\mathbb{R}^k$ such that $\langle\alpha,q\rangle+t\geq0$, we have
		\begin{description}
			\item[i.] $\mathcal{H}^{\langle\alpha,q\rangle+t+k\delta}({\underline{X}}_\alpha)\leq\,C\mathcal{H}_{\mu,\varphi}^{q,t}({\underline{X}}_\alpha).$
			\item[ii.] $\mathcal{P}^{\langle\alpha,q\rangle+t+k\delta}({\underline{X}}_\alpha)\leq\,C\mathcal{P}_{\mu,\varphi}^{q,t}({\underline{X}}_\alpha).$
		\end{description}
	\end{description}
	where $C=C(\alpha,q,k,\delta)>0$ is a generic constant.
\end{lemma}
\hskip-10pt{\it Proof.} {\bf 1.} {$\mathbf{i.}$} We prove the first part. For $m\in\mathbb{N}^*$, consider the set
$$
\overline{X}_m^\alpha=\{\,x\in\overline{X}^\alpha;\,\,\displaystyle\frac{\log(\mu_j(B(x,r)))}{\varphi(r)}\leq\alpha_j+\displaystyle\frac{\delta}{q_j};\,\,
0<r<\displaystyle\frac{1}{m},\;\;1\leq\,j\leq\,k\,\}.
$$
Let next $0<\eta<\displaystyle\frac{1}{m}$ and $(B(x_i,r_i))_i$ a centered $\eta$-covering of ${\overline X}_m^\alpha$. It holds that
$$
(\mu(B(x_i,r_i)))^{q}\geq\,e^{(\langle\alpha,q\rangle+k\delta)\varphi(r)}.
$$
Consequently, it holds from (\ref{condition-sur-phi}) that
$$
\mathcal{H}_\eta^{\langle\alpha,q\rangle+t+k\delta}(\overline{X}_m^\alpha)\leq\displaystyle\sum_i(2r_i)^{\langle\alpha,q\rangle+t+k\delta}
\leq\,C\displaystyle\sum_i(\mu(B(x_i,r_i)))^qe^{t\varphi(r)}.
$$
Hence, $\forall\eta>0$, there holds that
$$
\mathcal{H}_\eta^{\langle\alpha,q\rangle+t+k\delta}(\overline{X}_m^{\alpha})
\leq\,C\overline{\mathcal{H}}_{\mu,\varphi,\eta}^{q,t}(\overline{X}_m^{\alpha}).
$$
Which means that
$$
\mathcal{H}^{\langle\alpha,q\rangle+t+k\delta}(\overline{X}_m^{\alpha})
\leq\,C\overline{\mathcal{H}}_{\mu,\varphi}^{q,t}(\overline{X}_m^{\alpha})
\leq\,C\mathcal{H}_{\mu,\varphi}^{q,t}(\overline{X}_m^\alpha).
$$
Next, observing that $\overline{X}^\alpha=\displaystyle\bigcup_{m}\overline{X}_m^\alpha$, we obtain
$$
\mathcal{H}^{\langle\alpha,q\rangle+t+k\delta}(\overline{X}^\alpha)
\leq\,C\mathcal{H}_{\mu,\varphi}^{q,t}(\overline{X}^\alpha).
$$
{$\mathbf{ii.}$} For $q\in\mathbb{R}^{*,k}_+$ and $m\in\mathbb{N}^*$, consider the set $\overline{X}_m^\alpha$ defined previously and let $E\subseteq\overline{X}_m^\alpha$, $0<\eta<\displaystyle\frac{1}{m}$ and $\bigl(B(x_i,r_i)\bigr)_i$ a centered $\eta$-packing of $E$. We have
$$
\displaystyle\sum_i(2r_i)^{\langle\alpha,q\rangle+t+k\delta}\leq\,C\displaystyle\sum_i(\mu(B(x_i,r_i)))^qe^{t\varphi(r)}
\leq\,C\overline{\mathcal{P}}_{\mu,\varphi,\eta}^{q,t}(E).
$$
Consequently, $\forall\,\eta>0$,
$$
\overline{\mathcal{P}}_\eta^{\langle\alpha,q\rangle+t+k\delta}(E)\leq\,C\overline{\mathcal{P}}_{\mu,\varphi,\eta}^{q,t}(E).
$$
Hence, $\forall\,E\subseteq\overline{X}_m^\alpha$,
$$
\overline{\mathcal{P}}^{\langle\alpha,q\rangle+t+k\delta}(E)\leq\,C\overline{\mathcal{P}}_{\mu,\varphi}^{q,t}(E).
$$
Let next, $(E_i)_i$ be a covering of $\overline{X}_m^\alpha$. Thus,
$$
\begin{array}{lll}\mathcal{P}^{\langle\alpha,q\rangle+t+k\delta}(\overline{X}_m^\alpha)
&=&\mathcal{P}^{\langle\alpha,q\rangle+t+k\delta}\Biggl(\displaystyle\bigcup_i(\overline{X}_m^\alpha\cap\,E_i)\Biggr)\hfill\cr\medskip
&=&\displaystyle\sum_i\mathcal{P}^{\langle\alpha,q\rangle+t+k\delta}\Bigl(\overline{X}_m^\alpha\cap\,E_i\Bigr)\hfill\cr\medskip
&\leq&\displaystyle\sum_i\overline{\mathcal{P}}^{\langle\alpha,q\rangle+t+k\delta}\Bigl(\overline{X}_m^\alpha\cap\,E_i\Bigr)\hfill\cr\medskip
&\leq&C\displaystyle\sum_i\overline{\mathcal{P}}_{\mu,\varphi}^{q,t}\Bigl(\overline{X}_m^\alpha\cap\,E_i\Bigr)\hfill\cr\medskip
&\leq&C\displaystyle\sum_i\overline{\mathcal{P}}_{\mu,\varphi}^{q,t}(E_i).\end{array}
$$
Hence, $\forall,m$,
$$
\mathcal{P}^{\langle\alpha,q\rangle+t+k\delta}(\overline{X}_m^\alpha)\leq\,C\mathcal{P}_{\mu,\varphi}^{q,t}(\overline{X}_m^\alpha).
$$
Consequently,
$$
\mathcal{P}^{\langle\alpha,q\rangle+t+k\delta}(\overline{X}^\alpha)\leq\,C\mathcal{P}_{\mu,\varphi}^{q,t}(\overline{X}^\alpha).
$$
{\bf 2.} {$\mathbf{i.}$} and {$\mathbf{ii.}$} follow similar arguments and techniques as previously.
\begin{proposition}\label{spectreproposition1}
	Let $\alpha\in\mathbb{R}^k_+$, $q\in\mathbb{R}^k$ and $\varphi$ satisfying (\ref{condition-sur-phi}). The following assertions hold.
	\begin{description}
		\item[a.] Whenever $\langle\alpha,q\rangle+b_{\mu,\varphi}(q)\geq0$, we have
		\begin{description}
			\item[i.] $\hbox{dim}{\overline{X}}^\alpha\leq\langle\alpha,q\rangle+b_{\mu,\varphi}(q)$,\quad$\forall\,q\mathbb{R}^k_+$.
			\item[ii.] $\hbox{dim}{\underline{X}}_\alpha\leq\langle\alpha,q\rangle+b_{\mu,\varphi}(q)$,\quad$\forall\,q\mathbb{R}^k_-$.
		\end{description}
		\item[b.] Whenever $\langle\alpha,q\rangle+B_{\mu,\varphi}(q)\geq0$, we have
		\begin{description}
			\item[i.] $\hbox{Dim}{\overline{X}}^\alpha\leq\langle\alpha,q\rangle+B_{\mu,\varphi}(q)$,\quad$\forall\,q\mathbb{R}^k_+$.
			\item[ii.] $\hbox{Dim}{\underline{X}}_\alpha\leq\langle\alpha,q\rangle+B_{\mu,\varphi}(q)$,\quad$\forall\,q\mathbb{R}^k_-$.
		\end{description}
	\end{description}
\end{proposition}
\hskip-10pt{\bf Proof.} {\bf a. i.} It follows from Lemma \ref{spectrelemme1}, assertion {\bf 1. i.},
$$
\mathcal{H}^{\langle\alpha,q\rangle+t+k\delta}({\overline{X}}^\alpha)=0,\qquad\forall\,t>b_{\mu,\varphi}(q),\,\,\delta>0.
$$
Consequently,
$$
\hbox{dim}{\overline{X}}^\alpha\leq\langle\alpha,q\rangle+t+k\delta,\qquad\forall\,t>b_{\mu,\varphi}(q),\;\;\delta>0.
$$
Hence,
$$
\hbox{dim}\,\overline{X}^\alpha\leq\langle\alpha,q\rangle+b_{\mu,\varphi}(q).
$$
{\bf a. ii.} It follows from Lemma \ref{spectrelemme1}, assertion {\bf 2. i.}, as previously, that
$$
\mathcal{H}^{\langle\alpha,q\rangle+t+k\delta}(\underline{X}^\alpha)=0,\qquad\forall\,t>b_{\mu,\varphi}(q),\,\,\delta>0.
$$
Hence,
$$
\hbox{dim}{\underline{X}}_\alpha\leq\langle\alpha,q\rangle+t+k\delta,\qquad\forall\,t>b_{\mu,\varphi}(q),\;\;\delta>0
$$
and finally,
$$
\hbox{dim}{\underline{X}}_\alpha\leq\langle\alpha,q\rangle+b_{\mu,\varphi}(q).
$$
{\bf b. i.} observing Lemma \ref{spectrelemme1}, assertion {\bf 1. ii.}, we obtain
$$
\mathcal{P}^{\langle\alpha,q\rangle+t+k\delta}(\overline{X}^\alpha),\qquad\forall\,t>B_{\mu,\varphi}(q),\;\;\delta>0.
$$
Consequently,
$$
\hbox{Dim}\,\overline{X}^\alpha\leq\langle\alpha,q\rangle+t+k\delta,\qquad\forall\,t>B_{\mu,\varphi}(q),\;\;\delta>0.
$$
Hence,
$$
\hbox{Dim}\,\overline{X}^\alpha\leq\langle\alpha,q\rangle+B_{\mu,\varphi}(q).
$$
{\bf b. ii.} observing Lemma \ref{spectrelemme1}, assertion {\bf 2. ii.}, we obtain
$$
\mathcal{P}^{\langle\alpha,q\rangle+t+k\delta}(\underline{X}_\alpha)=0,\qquad\forall\,t>B_{\mu,\varphi}(q),\;\;\delta>0.
$$
Hence,
$$
\hbox{Dim}\,\underline{X}_\alpha\leq\langle\alpha,q\rangle+t+k\delta,\qquad\forall\,t>B_{\mu,\varphi}(q),\;\;\delta>0
$$
and finally,
$$
\hbox{Dim}\,\underline{X}_\alpha\leq\langle\alpha,q\rangle+B_{\mu,\varphi}(q).
$$
\begin{lemma}\label{spectrelemme2}
	$\forall\,q\in\mathbb{R}^k$ such that $\langle\alpha,q\rangle+b_{\mu,\varphi}(q)<0$ or $\langle\alpha,q\rangle+B_{\mu,\varphi}(q)<0$, we have $X(\alpha)=\emptyset$.
\end{lemma}
\hskip-10pt{\bf Proof.} It is based on
\begin{description}
	\item Claim 1. For $q\in\mathbb{R}^k_-$ with $\langle\alpha,q\rangle+b_{\mu,\varphi}(q)<0$ or $\langle\alpha,q\rangle+B_{\mu,\varphi}(q)<0$, $\underline{X}_\alpha=\emptyset$.
	\item Claim 2. For $q\in\mathbb{R}^k_+$ with $\langle\alpha,q\rangle+b_{\mu,\varphi}(q)<0$ or $\langle\alpha,q\rangle+B_{\mu,\varphi}(q)<0$, $\overline{X}^\alpha=\emptyset$.
\end{description}
Indeed, let $q\in\mathbb{R}^k_-$ and assume that $\underline{X}_\alpha\not=\emptyset$. This means that there exists at least one point $x\in\hbox{supp}(\mu)$ for which ${\underline\alpha}_{\mu_j}(x)\geq\alpha_j$, for $1\leq\,j\leq\,k$. Consequently, for all $\varepsilon>0$, there is a sequence $(r_n)_n\downarrow0$ and satisfying
$$
0<r_n<\displaystyle\frac{1}{n}\quad\hbox{and}\quad\mu_j(B(x,r_n))<e^{(\alpha_j-\varepsilon)\varphi(r_n)},\;\;1\leq\,j\leq\,k.
$$
Hence,
$$
\Bigl(\mu(B(x,r_n))\Bigr)^qe^{t\varphi(r_n)}>Ce^{(\langle(\alpha-\varepsilon\mathbb{I}),q\rangle+t)\varphi(r_n)}.
$$
Choosing $t=\langle(\varepsilon\mathbb{I}-\alpha),q\rangle$, this induces that $\mathcal{H}_{\mu,\varphi}^{q,t}(\{x\})>C>0$ and consequently,
$$
b_{\mu,\varphi}(q)\geq\,\hbox{dim}_{\mu,\varphi}^q(\{x\})\geq\,t,\quad\forall\,\varepsilon>0.
$$
Letting $\varepsilon\downarrow0$, it results that $b_{\mu,\varphi}(q)\geq-\langle\alpha,q\rangle$ which is impossible. So as the first part of Claim 1. The remaining part as well as Claim 2 can be checked by similar techniques.
\begin{theorem}\label{spectretheorem1}
	Let $\mu=(\mu_1,\,\mu_2,\,\dots,\,\mu_k)$ be a vector-valued Borel probability measure on $\mathbb{R}^d$ and $q\in\mathbb{R}^k$ fixed. Let further $t\in\mathbb{R}$, ${\underline{K}},\,{\overline{K}}>0$, $\nu$ a Borel probability measure supported by $\hbox{supp}(\mu)$, $\varphi:\,\mathbb{R}_+\rightarrow\mathbb{R}$ satisfying (\ref{condition-sur-phi}). Let finally $(r_{n})_n\subset]0,1[\downarrow0$ and satisfying
	$$
	\displaystyle\frac{\varphi(r_{n+1})}{\varphi(r_{n})}\rightarrow1\quad\hbox{and}\quad\displaystyle\sum_ne^{\varepsilon\varphi(r_{n})}<\infty,\
	\forall\varepsilon>0.
	$$
	Assume next the following assumptions.
	\begin{description}
		\item[A1.] $\forall\,x\in\hbox{supp}(\mu)$ and $r$ small enough,
		$$
		{\underline K}\leq\displaystyle\frac{\nu(B(x,r))}{\Bigl(\mu(B(x,r))\Bigr)^qe^{t\varphi(r)}}\leq{\overline K}.
		$$
		\item[A.2] $C(p)=\displaystyle\lim_{n\rightarrow+\infty}C_{n}(p)$ exists and finite for all $p\in\mathbb{R}$, where
		$$
		C_{n}(p)=\displaystyle\frac{-1}{\varphi(r_{n})}\log\biggl(\displaystyle\int_{supp(\mu)}\Bigl(\mu(B(x,r_{n}))\Bigr)^pd\nu(x)\biggr).
		$$
	\end{description}
	Denote next $\alpha_-^0=-\nabla_-C(0)$, $\alpha_+^0=-\nabla_+C(0)$ and $\Psi_q(a,b)=aq+b$, $\forall\,a,b$. The following assertions hold.
	\begin{description}
		\item[i.] For $q\in\mathbb{R}^k_-$, we have
		$$
		\hbox{dim}(\underline{X}_{\alpha_+^0}\cap\overline{X}^{\alpha_-^0})
		\geq\Psi_q(\alpha_-^0,\Lambda_\mu(q))\geq\Psi_q(\alpha_-^0,B_\mu(q))\geq\Psi_q(\alpha_-^0,b_\mu(q)).
		$$
		For $q\in\mathbb{R}^k_+$,
		$$
		\hbox{dim}(\underline{X}_{\alpha_+^0}\cap\overline{X}^{\alpha_-^0})
		\geq\Psi_q(\alpha_+^0,\Lambda_\mu(q))\geq\Psi_q(\alpha_+^0,B_\mu(q))\geq\Psi_q(\alpha_+^0,b_\mu(q)).
		$$
		\item[ii.] Whenever $C$ is differentiable at 0, we have
		$$
		f_\mu(-\nabla\,C(0))=b_\mu^*(-\nabla\,C(0))=B_\mu^*(-\nabla\,C(0))=\Lambda_\mu^*(-\nabla\,C(0)).
		$$
	\end{description}
\end{theorem}
\hskip-10ptThe proof of this result is based on the application of a large deviation formalism. This will permit to obtain a measure $\nu$ supported by ${\underline X}_{-\nabla_+C(0)}\cap{\overline X}^{-\nabla_-C(0)}$. To do this, we re-formulate a mixed large deviation formalism to be adapted to the mixed multifractal formalism raised in our work.\\
{\bf Proof of Theorem \ref{spectretheorem1}.} For $x\in supp(\mu)$, let
$$
{\underline\alpha}_{\mu_j}^{\varphi}(x,r_n)=\displaystyle\liminf_n\displaystyle\frac{\log\Bigl[\mu_j(B(x,r_n)\Bigr]}{\varphi(r_n)}
$$
and
$$
{\overline\alpha}_{\mu_j}^{\varphi}(x,r_n)=\displaystyle\limsup_n\displaystyle\frac{\log\Bigl[\mu_j(B(x,r_n)\Bigr]}{\varphi(r_n)}.
$$
{\bf i.} Using the hypothesis {\bf A1.} and Lemma \ref{anouarlemme3} we obtain
$$
b_{\mu,\varphi}(q)=B_{\mu,\varphi}(q)=\Lambda_{\mu,\varphi}(q)=t.
$$
Next, it is straightforward that the set
$$
M=\left\{\ x\in supp(\mu)\ ;\ -\nabla_+C(0)\leq{\underline{\alpha}}_{\mu,\varphi}(x,r_n)\leq{\overline{\alpha}}_{\mu,\varphi}(x,r_n)\leq-\nabla_-C(0)\ \right\}
$$
coincides with $\displaystyle{\underline X}_{-\nabla_+C(0)}\cap{\overline X}^{-\nabla_-C(0)}$. Hence, by setting in the mixed large deviation formalism \ref{largedeviation},
$$
\Omega=supp(\mu),\;\;{\mathcal A}={\mathcal B}(supp(\mu)),
$$
$$
I\!\!P=\mu,\;\;W_n(x)=\log(\mu(B(x,r_n)))
$$
and
$$
a_n=-\varphi(r_n),
$$
it holds that
$$
{\underline\alpha}_{\mu,\varphi}(x)\geq\left\{\begin{array}{lll}-\nabla_-C(0)q+t&\hbox{for}&q\leq0\hfill\cr\medskip
-\nabla_+C(0)q+t&\hbox{for}&q\geq0.\end{array}\right.
$$
Finally, applying the famous Billingsley's Theorem \cite{Billingsley}, we obtain
$$
\hbox{dim}\,M\geq\left\{\begin{array}{lll}-\nabla_-C(0)q+t&\hbox{for}&q\leq0\hfill\cr\medskip
-\nabla_+C(0)q+t&\hbox{for}&q\geq0.\end{array}\right.
$$
{\bf ii.} Remark that if $C$ is differentiable at 0, item {\bf i.} states that
$$
\hbox{dim}\,M\geq-\nabla\,C(0)q+t\geq\Lambda_\mu^*(-\nabla\,C(0))\geq B_\mu^*(-\nabla\,C(0))\geq b_\mu^*(-\nabla\,C(0)).
$$
In the other hand, since the set $M$ is not empty, Lemma \ref{spectrelemme2} implies that
$$
-\nabla\,C(0)q+t\geq0.
$$
Hence, Proposition \ref{spectreproposition1} yields that
$$
\hbox{dim}\,M\leq-\nabla\,C(0)q+t
$$
for any $q\in\mathbb{R}^k$. Thus, taking the inf on $q$, we obtain
$$
\hbox{dim}\,M\leq b_\mu^*(-\nabla\,C(0))\leq\,B_\mu^*(-\nabla\,C(0))\leq\Lambda_\mu^*(-\nabla\,C(0)).
$$
{\bf iii.} We firstly claim that, there exists $\beta>0$ such that, for all $x\in supp(\mu)$ and $0<r<<<1$, we have
$$
\displaystyle\frac{\mu(B(x,2r))}{\mu(B(x,r))}<\beta.
$$
So let $(B(x_{ij},r_n))_{1\leq i\xi, j}$ the $\xi$ sets relatively to Besicovitch theorem extracted from the set $(B(x_i,r_n))_i$. A careful computation yields that
\begin{equation}\label{equation1spectre}
	|p+q-\mathbb{I}|I_\mu^{p+q-\mathbb{I}}=C_q(p)+t_q;\quad\forall\,p,q\in\mathbb{R}^k
\end{equation}
where $|p+q-\mathbb{I}|=\displaystyle\sum_{i=1}^k(p_i+q_i-1)$. Theorem \ref{anouartheorem1} and Proposition \ref{anouarproposition5} guarantees that
$$
|p+q-\mathbb{I}|I_\mu^{p+q-1}=C_\mu^{p+q}(supp(\mu))=\Lambda_\mu(p+q).
$$
Consequently,
$$
C_q(p)=\Lambda_\mu(p+q)-\Lambda_\mu(p).
$$
So, if $\Lambda_\mu$ is differentiable at $q$, $C_q$ will be too at 0 and $\nabla\,C_q(0)=\nabla\Lambda_\mu(q)$. Thus, using the mixed large deviation formalism, we obtain
$$
\alpha_\mu(x)=-\nabla\,C_q(0)\ ;\ \nu_q\ \hbox{for almost all}\ x\in
supp(\mu).
$$
hence, finally, $\alpha_\mu(x)=-\nabla\Lambda_\mu(q)$.\\
{\bf iv.} Let $q$ be such that $\nabla\Lambda_\mu(q)$ exists. Then $\nabla\,C_q(0)$ exists too. So, item {\bf ii.} states that
$$
f_\mu(-\nabla\,C(0))=\Lambda_\mu^*(-\nabla\,C(0)).
$$
Which completes the proof.
\begin{theorem}\label{spectretheorem2}
	Assume that the hypotheses of Theorem \ref{spectretheorem1} are satisfied for all $q\in\mathbb{R}^k$. Then, the following assertions hold.
	\begin{description}
		\item[i.] $\alpha_\mu=-B_\mu,\quad\nu_q\,\,a.s$, whenever $B_\mu$ is differentiable at $q$.
		\item[ii.] $Dom(B)\subseteq\alpha_\mu(supp(\mu))$ and $f_\mu=B_\mu^*$ on $Dom(B)$.
	\end{description}
\end{theorem}
\hskip-10pt{\bf Proof.}\\
{\bf i.} Using the same notations as in Theorem \ref{spectretheorem1}, we obtain $C$ differentiable at 0, $B_{\mu,\varphi}$ differentiable at $q$, and $\nabla\,C(0)=\nabla\,B_{\mu,\varphi}(q)$. In the other hand, we obtain also
$$
\alpha_\mu^\varphi(x)=\alpha_\mu^\varphi(x,r_n)=\displaystyle\lim_n\frac{W_n(x)}{-a_n}=-\nabla\,C(0)=\nabla\,B_\mu^\varphi(q),\quad\nu\;\;a.s.
$$
{\bf ii.} follows immediately from {\bf i.} and Theorem \ref{spectretheorem1}.
\section{Appendix}
\subsection{Besicovitch covering theorem}
\begin{theorem} \label{besicovitch}
	There exists a constant $\xi\in\mathbb{N}$ satisfying: For any $E\in\mathbb{R}^d$ and $(r_x)_{x\in\,E}$ a bounded set of positive real numbers, there exists $\xi$ sets $B_1$, $B_2$, ..., $B_{\xi}$, that are finite or countable composed of balls $B(x,r_x)$, $x\in\,E$ such that
	\begin{itemize}
		\item $E\subseteq\displaystyle\bigcup_{1\leq\,i\leq\xi}\displaystyle\bigcup_{B\in\,B_i}B$.
		\item each $B_i$ is composed of disjoint balls.
	\end{itemize}
\end{theorem}
\subsection{A mixed large deviation theorem}
To do this, we re-formulate a mixed large deviation formalism to be adapted to the mixed multifractal formalism raised in our work.
\begin{theorem}\label{largedeviation}
	Consider a sequence $(W_n=(W_{n,1},\,W_{n,2},\,\dots,\,W_{n,k}))_n$ of vector-valued random variables on a probability space $(\Omega,\,\mathcal{A},\,\mathbb{P})$ and $(a_n)_n\subset]0,+\infty[$ with $\displaystyle\lim_{n\rightarrow+\infty}a_n=+\infty$. Let next the function
	$$
	\begin{array}{lll}C_n&:&\mathbb{R}^k\rightarrow\overline{\mathbb{R}}\hfill\cr\medskip
	& &t\mapsto C_n(t)=\displaystyle\frac{1}{a_n}\log\Bigl(E(\exp(\langle\,t,W_n\rangle))\Bigr).\end{array}
	$$
	Assume that
	\begin{description}
		\item[A1.] $C_n(t)$ is finite for all $n$ and $t$.
		\item[A2.] $C(t)=\displaystyle\lim_{n\rightarrow+\infty}C_n(t)$ exists and is finite for all $t$.
	\end{description}
	There holds that
	\begin{description}
		\item[i.] The function $C$ is convex.
		\item[ii.] If $\nabla_-C(t)\leq\nabla_+C(t)<\alpha$, for some $t\in\mathbb{R}^k$, then
		$$
		\displaystyle\limsup_{n\rightarrow+\infty}\displaystyle\frac{1}{a_n}
		\log\Biggl(e^{-a_nC(t)}E\biggl(\exp(\langle\,t,W_n\rangle)1_{\{\frac{W_n}{a_n}\geq\alpha\}}\biggr)\Biggr)<0.
		$$
		\item[iii.] If $\displaystyle\sum_ne^{-\varepsilon a_n}<\infty$ for all $\varepsilon>0$, then
		$$
		\displaystyle\limsup_{n\rightarrow+\infty}\frac{W_n}{a_n}\leq\nabla_+C(0)\qquad \mathbb{P}\ a.s.
		$$
		\item[iv.] If $\alpha<\nabla_-C(t)\leq\nabla_+C(t)$, for some $t\in\mathbb{R}^k$, then
		$$
		\displaystyle\limsup_{n\rightarrow+\infty}\displaystyle\frac{1}{a_n}
		\log\Biggl(e^{-a_nC(t)}E\biggl(\exp(\langle\,t,W_n\rangle)1_{\{\frac{W_n}{a_n}\leq\alpha\}}\biggr)\Biggr)<0.
		$$
		\item[v.] If $\displaystyle\sum_ne^{-\varepsilon a_n}$ is finite for all $\varepsilon>0$, then
		$$
		\nabla_-C(0)\leq\displaystyle\limsup_{n\rightarrow+\infty}\frac{W_n}{a_n}\qquad\mathbb{P}\ a.s.
		$$
	\end{description}
\end{theorem}
\hskip-10pt{\bf Proof.}\\
{\bf i.} It follows from Holder's inequality.\\
{\bf ii.} Let $h\in\mathbb{R}^{*,k}_+$ be such that $C(t)+\langle\alpha,h\rangle-C(t+h)>0$. We have
$$
\begin{array}{lll}&&\displaystyle\frac{1}{a_n}\log\biggl[e^{-a_nC(t)}
\mathbb{E}\Bigl(\exp(\langle\,t,W_n\rangle)1_{\{\frac{W_n}{a_n}\geq\alpha\}}\Bigr)\biggr]\hfill\cr\medskip
&=&\displaystyle\frac{1}{a_n}\log\biggl[e^{-a_nC(t)}\displaystyle\int_{\{\frac{W_n}{a_n}\geq\alpha\}}
e^{\langle\,t,W_n\rangle}d\mathbb{P}\biggr]\hfill\cr\medskip
&=&\displaystyle\frac{1}{a_n}\log\biggl[e^{-a_n(C(t)+\langle\alpha,h\rangle)}\displaystyle\int_{\{\frac{W_n}{a_n}\geq\alpha\}}
e^{\langle\,t,W_n\rangle+a_n\langle\alpha,h\rangle}d\mathbb{P}\biggr]\hfill\cr\medskip
&\leq&\displaystyle\frac{1}{a_n}\log\biggl[e^{-a_n(C(t)+\langle\alpha,h\rangle)}\displaystyle\int_{\{\frac{W_n}{a_n}\geq\alpha\}}
e^{\langle\,t+h,W_n\rangle}d\mathbb{P}\biggr]\hfill\cr\medskip
&\leq&\displaystyle\frac{1}{a_n}\log\biggl[e^{-a_n(C(t)+\langle\alpha,h\rangle)}\mathbb{E}
(\exp(\langle\,t+h,W_n\rangle))\biggr]\hfill\cr\medskip
&=&\displaystyle\frac{1}{a_n}\log\biggl[e^{-a_n(C(t)+\langle\alpha,h\rangle-C_n(t+h))}\biggr]\hfill\cr\medskip
&=&-(C(t)+\langle\alpha,h\rangle-C_n(t+h)).\end{array}
$$
Next, by taking the limsup as $n\longrightarrow+\infty$, the result follows immediately.\\
{\bf iii.} Denote for $n,m\in\mathbb{N}$,
$$
T_{n,m}=\{\frac{W_n}{a_n}\geq\nabla_+C(0)+\frac{1}{m}\}.
$$
By choosing in item {\bf ii.} $t=0$ and $\alpha=\nabla_+C(0)+\frac{1}{m}$, and observing that $C(0)=0$, we obtain
$$
\displaystyle\limsup_{n\rightarrow+\infty}\displaystyle\frac{1}{a_n}
\log\Biggl(E\biggl(1_{\{\frac{W_n}{a_n}\geq\nabla_+C(0)+\frac{1}{m}\}}\biggr)\Biggr)<0.
$$
which means that
$$
\displaystyle\limsup_n\frac{1}{a_n}\log\mathbb{P}(T_{n,m})<0.
$$
Consequently, for some $\varepsilon>0$ and $n$ large enough, there holds that
$$
\displaystyle\limsup_n\frac{1}{a_n}\log\mathbb{P}(T_{n,m})<-\varepsilon.
$$
Thus,
$$
\mathbb{P}(T_{n,m})<e^{-\varepsilon a_n}
$$
which implies the convergence of the series $\displaystyle\sum_n\mathbb{P}(T_{n,m})$. Hence, using Borel-Cantelli theorem, we obtain $$
\mathbb{P}(\displaystyle\limsup_nT_{n,m})=0,\;\forall m.
$$
Therefore,
$$
\mathbb{P}\biggl(\displaystyle\limsup_n\displaystyle\frac{W_n}{a_n}>\nabla_+C(0)\biggr)=\mathbb{P}(\displaystyle\bigcup_m\limsup_nT_{n,m})=0
$$
and finally,
$$
\displaystyle\limsup_n\displaystyle\frac{W_n}{a_n}\leq\nabla_+C(0),\;\mathbb{P}\,\hbox{a.s}.
$$
\section*{Acknowledgment} The second author would like to thank Professor Lars Olsen for the first reading of this work and for the interest he gave to it.

\end{document}